\pdfoutput=1
\documentclass[9pt,twoside,epsf]{article}
\usepackage{stmaryrd}
\usepackage{amsmath,latexsym,amssymb,amsfonts,amsbsy}
\usepackage{bm}
\usepackage{graphicx,subfigure}
\usepackage{cases}
\usepackage{multirow}
\usepackage{booktabs}
\usepackage{diagbox}
\usepackage{float}
\usepackage{enumitem}
\usepackage[hidelinks]{hyperref}
\usepackage{cleveref}

\usepackage{lineno}
\usepackage{xcolor}
\newfont{\bb}{msbm10}
\def\Bbb#1{\mbox{\bb #1}}

\def\T{\top}

\def\diag{{\rm diag}}

\def\E{\mathbb{E}}
\def\rank{{\rm rank}}
\def\range{{\rm range}}

\def\zspace{{\rm null}}

\def\pr{{\rm Pr}}

\def\cov{{\rm Cov}}

\usepackage{colortbl}
\usepackage{ragged2e}
\usepackage[margin=2em,labelsep=space,skip=0.5em,font=normalfont]{caption}
\usepackage{algorithm}
\usepackage{algpseudocode}
\floatname{algorithm}{\color{black} Algorithm}
\algrenewcommand{\algorithmiccomment}[1]{\quad{\color{red}\%\ #1}}
\numberwithin{algorithm}{section}
\makeatletter
\newenvironment{breakalgo}[2]{%
  \captionsetup{margin=0pt,justification=RaggedRight,singlelinecheck=false}%
  \par\noindent%
  \medskip%
  \rule{\linewidth}{0.8pt}%
  \vspace{-0.5\baselineskip}%
  \noindent\captionof{algorithm}{#1}\label{#2}%
  \vspace{-0.7\baselineskip}%
  \noindent\rule{\linewidth}{.4pt}%
  \vspace{-0.3\baselineskip}%
}{%
  \vspace{-.75\baselineskip}%
  \rule{\linewidth}{.4pt}%
  \medskip%
}

\makeatother

\newtheorem{example}{Example}[section]

\newtheorem{remark}{Remark}[section]
\newtheorem{theorem}{Theorem}[section]
\newtheorem{definition}{Definition}[section]
\newtheorem{lemma}{Lemma}[section]
\newtheorem{corollary}{Corollary}[section]

\newcommand{\reals}{\makebox{{\Bbb R}}}

\newcommand{\bbb}[1]{\text{\bf #1}}
\newcommand{\txt}[1]{\text{\rm #1}}
\newcommand{\txttiny}[1]{\text{\tiny\rm #1}}

\makeatletter 

\@addtoreset{equation}{section}
\makeatother 

\begin{document}
\cleardoublepage
\pagestyle{myheadings}

\bibliographystyle{plain}

\title{\bf{Randomized batch-sampling Kaczmarz methods for solving linear systems}
\thanks{First author: Dong-Yue Xie (xiedongyue@nuaa.edu.cn); Corresponding author: Xi Yang (yangxi@nuaa.edu.cn).}}

\author{Dong-Yue Xie\thanks{School of Mathematics, Nanjing University of Aeronautics and Astronautics, Nanjing 211106, China.}, Xi Yang\footnotemark[2]}
	
\maketitle

\begin{abstract}
To conduct a more in-depth investigation of randomized solvers for solving linear systems, we adopt a unified randomized batch-sampling Kaczmarz framework with per-iteration costs as low as cyclic block methods, and develop a general analysis technique to establish its convergence guarantee. With concentration inequalities, we derive new expected linear convergence rate bounds. The analysis applies to any randomized non-extended block Kaczmarz methods with arbitrary static stochastic samplings. In addition, the new rate bounds are scale-invariant, which eliminate the dependence on the magnitude of the data matrix. In most experiments, the new bounds are significantly tighter than existing ones and better reflect the empirical convergence behavior of block methods. Within this new framework, the batch-sampling distribution, as a learnable parameter, provides the possibility for block methods to achieve efficient performance in specific application scenarios, which deserves further investigation.

\bigskip



\end{abstract}

\section{Introduction}
\label{sec:intro}

Consider large-scale linear systems of the form
\begin{eqnarray}
\label{linear-system}
  Ax &=& b,\quad\mbox{with $A\in\reals^{m\times n}$, $b\in\reals^m$, and $x\in\reals^n$},
\end{eqnarray}
where the data matrix $A$ in (\ref{linear-system}) appears in many real-world applications. Such linear systems have always been a fundamental and important scientific issue in the fields of science and industry, including computer tomography \cite{Hansen2021CT}, partial differential equations \cite{Olshanskii2014Iterative}, image reconstruction \cite{Byrne2003Unified}, signal processing \cite{PopaZdunek04} and machine learning \cite{Chang2008CD,Patrascu2018SPP}. Since large-scale linear systems are often too massive to be handled efficiently by direct methods, researchers have been motivated to develop iterative algorithms. Among various approaches, the Kaczmarz method \cite{Kaczmarz37}, originally proposed in 1937, has received significant attention as a typical row-action method due to its low per-iteration cost and ease of implementation. At each iteration, the simple Kaczmarz projects the current iterate onto a hyperplane defined by a single row of the system. Specifically, the update is defined by
\begin{eqnarray}
\label{RK-method}
  x^{(k+1)} &=& x^{(k)} + \frac{(b_{j_k}-A^{(j_k)}x^{(k)})}{\|A^{(j_k)}\|_2^2}\left(A^{(j_k)}\right)^{\T},
\end{eqnarray}
where $A^{(j)}$ denotes the $j$th row of $A$, $b_j$ denotes the $j$th entry of $b$, and $x^{(k)}$ denotes the $k$th iterate. Its randomized variant, known as the randomized Kaczmarz (RK) method \cite{Strohmer09}, is particularly appealing for large-data problems, where randomization leads to provable linear convergence guarantees. Recent publications on randomized versions of the Kaczmarz method have been triggered by the work of Strohmer and Vershynin \cite{Strohmer09}, which established the foundation for subsequent variants \cite{BaiWu18SIAM,AML-BaiWu18,BaiWang2023} and extensions \cite{LAA-BaiWu19,BaiWu21SIAM,Zouzias13,Ma15} of randomized Kaczmarz-type algorithms.

The block Kaczmarz method further enhances the computational efficiency of the simple RK method by simultaneously projecting the iterate onto the solution spaces of multiple equations. Specifically, we denote by $\tau \subset \{1,\dots,m\}$ a subset of row indices, $A_\tau$ the sub-data block of $A$ with rows indexed by $\tau$, and $b_\tau$ the corresponding part of the vector $b$. Starting from an initial guess $x^{(0)}$, if an index subset $\tau^{(k)}$ is cyclically or randomly sampled at iteration $k$, the next iterate is then obtained by projecting the current iterate onto the solution space $\{x\in\mathbb{R}^n:A_{\tau^{(k)}}x=b_{\tau^{(k)}}\}$, which leads to the following update
\begin{eqnarray}
	\label{BRK-method}
	x^{(k+1)} &=& x^{(k)} + A_{\tau^{(k)}}^{\dag}r_{\tau^{(k)}}^{(k)}, \ \text{with $r_{\tau^{(k)}}^{(k)}=b_{\tau^{(k)}}-A_{\tau^{(k)}}x^{(k)}$,}
\end{eqnarray}
where $A_{\tau^{(k)}}^{\dag}$ denotes the pseudoinverse of $A_{\tau^{(k)}}$. The block methods are advantageous when the runtime to process a block of equations jointly in one iteration (highly suitable for parallel computing in computational architectures) is significantly less than performing separate single row updates \cite{Dekel2012Optimal,ShalevShwartzZhang2013}.

Bai and Liu \cite{Bai13} proved the convergence of the cyclic block Kaczmarz (CBK) method via the Meany inequality, while Needell and Tropp \cite{NeedellTropp2014LAA} obtained the expected linear convergence rate bound of the randomized block Kaczmarz (RBK) method for row-normalized data matrices (a refined convergence guarantee for RBK
without the requirement of row-normalization is provided by Lok and Rebrova in \cite{LokRebrova2025X}). The RBK method has inspired a variety of accelerated or extended randomized block Kaczmarz methods \cite{gower2015randomized,Necoara19SIMAX,Chen2022DBK,DuSiSun20SISC,Tan2025ADBK}. More recently, Gower et al. \cite{gower2021adaptive} studied the sketch-and-project (SAP) framework, which is built upon randomized sketching techniques \cite{Woodruff2014,MartinssonTropp2020}. Within the SAP framework, both the RK and RBK methods can be viewed as special cases. Gower et al. also proposed adaptive block variants under the SAP framework and provided corresponding convergence analysis, related work can be found in \cite{gower2015randomized,gower2021adaptive,DerezinskiRebrova2024}. In these works, greedy and adaptive variants are shown to attain faster theoretical convergence rates, but determining the projection subspace at each iteration, that is, selecting the block indices, can be computationally expensive. In practice, this step may require global scans or scattered memory access that violate data locality \cite{Denning2005}, since evaluating all block residuals requires traversing multiple chunks of data, which causes non-contiguous memory access and extra data movement. As a result, the actual performance of these algorithms depends not only on the theoretical convergence rate but also on the per-iteration cost. Based on the above considerations, this paper focuses on RBK methods with non-adaptive sampling rules, while the developed analysis technique may be generalized to greedy and adaptive variants, with related work left for future research.

In addition to the development of new RK-type methods, the convergence analysis of such methods is also a significant research direction. Existing convergence analyses of RK-type methods are usually established under worst-case assumptions. As a result, the worst-case bounds often provide only conservative estimates of the actual convergence behavior. For these reasons, some further theoretical research has pursued more refined analyses of the RK-type methods through replacing worst-case contraction by typical-case behavior or empirical averages. Chen and Powell \cite{ChenPowell2012JFAA}, along with Lin and Zhou \cite{LinZhou2015JMLR} provided considerably stronger (almost sure) convergence guarantees for the RK method under stronger distributional assumptions.  For the tall coefficient matrix $A$, Wang et al. \cite{WangAgaskarLu2015SampTA} derived an exact mean squared error (MSE) analysis for the RK method. Building on this exact MSE viewpoint, Bai and Wu \cite{LAA-BaiWu18} established exact formulas for the MSE of both tall and wide matrices, and further derived sharper upper bounds on the convergence rate than the classical estimate. Anderson et al. \cite{Anderson2025Arxiv}  provided confidence intervals around expected error bounds for the RK method, and further established high-probability results for the error along the whole random trajectory, where the analysis is based on upper bounds for the variance and concentration of the error. However, these refined analytical results are all concerned with the simple RK methods rather than its block variants.

\subsection{Motivation and contributions}
The present work is mainly motivated by the following two observations. Firstly, existing block Kaczmarz methods generally describe block sampling rules at a macroscopic level. Such rules are manually specified and are generally non-optimal, which restricts the attainable performance of randomized methods. For instance, one typically first partitions the matrix rows into fixed blocks and then defines a probability distribution over these blocks, or selects a working block at each iteration according to some prescribed sampling criterion. While such formulations are effective for particular sampling mechanisms, they are not sufficiently flexible to support a unified characterization of general stochastic selection rules. Secondly, randomized linear solvers often exhibit a substantial gap between theoretical convergence rate bounds and empirical convergence rates. In particular, existing expected convergence analyses for block Kaczmarz methods are almost always established under worst-case assumptions. Consequently, the resulting theoretical convergence rates are often overly conservative and may not accurately capture the actual convergence performance of the algorithms in practice.

Motivated by these observations, the main contributions and novel perspectives of this paper are summarized as follows.
\begin{itemize}[itemsep=2pt, topsep=0pt, parsep=0pt, partopsep=0pt]
	\item[\textbf{1.}] \textbf{A unified RBSK framework for any static stochastic batch sampling.}
	
	We propose a unified RBSK framework that parameterizes static stochastic batch-sampling rules by a joint distribution $\bbb{P}$, with the detailed formulation and corresponding algorithm presented in Section \ref{sec:RBSK} (Algorithm \ref{alg-RBSK}). The core idea is to describe a batch-sampling rule through a joint distribution $\bbb{P}$ over row batches, where the batch selected at each iteration is represented by a random vector of row indices. To the best of our knowledge, this is the first work to characterize sampling rules in this manner. Owing to this fine-grained definition of sampling distributions, the RBSK framework is sufficiently general to characterize any randomized non-extended block Kaczmarz method employing arbitrary static stochastic samplings, while maintaining the same low per-iteration cost as the CBK and RBK methods.
	
	\item[\textbf{2.}] \textbf{A unified convergence theory with sharper and scale-invariant bounds.}
	
	Based on this framework, we establish a unified convergence theory and derive new convergence rate bounds. In particular, our results include worst-case bounds and an improved non-worst-case bound obtained via concentration inequalities. As an important part of our analysis, we introduce a scaling operator $S$ and obtain scale-invariant bounds, so that the resulting convergence estimates no longer depend on the magnitude of the data matrix. The discussion in Remarks \ref{remark5.1} and \ref{remark5.2} further clarifies the distinct roles of $\bbb{P}$ and $S$:
\begin{itemize}
  \item the joint distribution $\bbb{P}$ determines the actual performance of the randomized methods;
  \item the scaling operator $S$ only affects the tightness of convergence rate bounds in theory.
\end{itemize}
Moreover, under the row-paving setting (suggested by Needell and Tropp \cite{NeedellTropp2014LAA}), a classical bound for the RBK method reads that
\begin{align*}
  1 - \lambda_{\min}\left(A^{\top} \frac{1}{\beta L} A\right),
\end{align*}
where $\beta$ denotes a paving-related quantity, while $L$ is the number of paved blocks. We establish the following new bound in this setting
\begin{align*}
  \min_{S}(1-\E\left(\xi_{\tau} \right)),\  \xi_{\tau}=\lambda_{\min}(W_{\tau}^{\T}A^{\T} D^2AW_{\tau}),\ \tau\sim\bbb{P},
\end{align*}
where $D$ (associated with $S$) and $W_{\tau}$ will be introduced in subsequent sections. The new bound features the following improvements: $D$ constitutes a refinement of the quantity $1/(\beta L)$; $W_{\tau}$ restricts the matrix to a subspace of $\reals^n$; taking the expectation over $\xi_{\tau}$ yields a non-worst-case bound; minimization over $S$ leads to a tighter bound. Owing to these features, the proposed bound is substantially tighter than the one established in \cite{NeedellTropp2014LAA}. In most numerical experiments on synthetic multi-scale and ill-conditioned systems as well as on sparse matrices from the SuiteSparse Collection \cite{davis2011university}, the new bound is validated to be tighter and more consistent with empirical convergence rates than the existing bounds from \cite{NeedellTropp2014LAA,gower2021adaptive}. The test bounds are listed in Table~\ref{tab:formulas}.
	
	\item[\textbf{3.}] \textbf{A novel optimization perspective via learnable sampling distributions.}
	
	A significant advantage of the RBSK framework is that it parameterizes the sampling rules, which offers the potential for randomized methods to approach their intrinsic performance limits. Specifically, the joint distribution $\bbb{P}$ can be treated as an optimizable and learnable parameter, which suggests a new way to construct efficient Kaczmarz-type methods through optimizing the sampling distribution, and provides the possibility to achieve highly efficient performance in specific application scenarios. A systematic study of this direction is beyond the scope of the present paper, which deserves further investigation.
\end{itemize}
\subsection{Organization}
The organization of this paper is as follows. In Section \ref{sec:prelimNota}, we introduce some necessary definitions, theorems, and notations. In Section \ref{sec:RBSK}, we present the detailed formulation of the RBSK method. In Section \ref{sec:lemma}, we establish several main lemmas that are useful for analyzing the convergence rate of the RBSK method. In Section \ref{sec:convAnal}, we present the new convergence analysis of the RBSK method. Numerical results are reported in Section \ref{sec:expResult}. Finally, Section \ref{sec:concl} ends the paper with summaries and directions for future research. Additional auxiliary lemmas and representative instances of stochastic batch-sampling rules are provided in the appendices.

\section{Preliminaries and Notations}
\label{sec:prelimNota}

This section introduces fundamental statistical quantities and concentration inequalities involved in convergence analysis, which are necessary tools for proving our main results. Additionally, this section provides a summary of notations used throughout the paper.

\subsection{Basic sample statistics}

\begin{definition}\label{def:basic-statistics}
	Let $X$ be a real-valued random variable defined on a probability space $(\Omega,\mathcal{F},\bbb{P})$, with finite mean $\mu=\mathbb{E}[X]$ and finite variance $\sigma^2=\operatorname{Var}(X)$. The distribution of $X$ is referred to as the population distribution, and $\mu, \sigma^2$ are called the population mean and population variance, respectively. A sample of size $n$ is a sequence of independent and identically distributed (i.i.d.) instances $\{x_i\}_{i=1}^n$ of $X$, referred to as a sample sequence. Based on the definitions of population and sample above, we introduce the following statistical quantities.
	\begin{itemize}
		
		\item \textbf{Sample Mean:}
		Given a sample sequence $\{x_i\}_{i=1}^n$, the sample mean is defined by
		\begin{eqnarray*}
			\bar{x}=\frac{1}{n} \sum_{i=1}^n x_i .
		\end{eqnarray*}
		It is an unbiased estimator of the population mean, i.e., $\mathbb{E}[\bar{x}] = \mu$.	
		
		\item \textbf{Sample Covariance:}
		Let $\{x_i\}_{i=1}^n$ and $\{y_i\}_{i=1}^n$ be two sample sequences of size $n$, with respective sample means $\bar{x}$ and $\bar{y}$. The sample covariance is defined by
		\begin{eqnarray*}
			\cov\left(x, y\right)=\frac{1}{n-1} \sum_{i=1}^n\left(x_i-\bar{x}\right)\left(y_i-\bar{y}\right).
		\end{eqnarray*}
		When the two sample sequences coincide, i.e., $x_i = y_i$ for all $i$, $\cov\left(x, y\right)$ reduces to the \textbf{Sample Variance} $s^2=\cov\left(x, x\right)$, and $\mathbb{E}[s^2] = \sigma^2$.
		
	\end{itemize}
\end{definition}

\subsection{Concentration inequalities}
\begin{theorem}[Hoeffding's inequality]
	\label{Hoeffding}
	Let $X_1, \ldots, X_n$ be independent random variables such that $X_i$ takes its value in $\left[a_i,b_i\right]$ almost surely for all $i \leq n$. Consider the sum of these random variables,
	\begin{eqnarray*}
		S_n=X_1+\cdots+X_n,
	\end{eqnarray*}
	then, for $\epsilon>0$, it holds that
	\begin{eqnarray*}
		\pr\left(\left|S_n-\E(S_n)\right|\geq\epsilon\right) &\leq& 2 \exp \left(-\frac{2 \epsilon^2}{\sum_{i=1}^n\left(b_i-a_i\right)^2}\right).
	\end{eqnarray*}
\end{theorem}
The inequality above is a classical result, and a proof can be found in probability and statistics-related references such as \cite{Vershynin2018,Boucheron2013}. A direct consequence of Theorem \ref{Hoeffding} for the sample mean is as follows.
\begin{corollary}
	\label{Hoeffding-co}
	Let $x_1, \dots, x_n$ be independent samples drawn from a distribution of a random variable $X$ supported on $[a, b]$ almost surely for all $i \leq n$. Let $\bar{x} = \frac{1}{n} \sum_{i=1}^n x_i$ be the sample mean. For $\varepsilon > 0$ and $\delta \in (0,1)$, if the number of samples satisfies
	\begin{eqnarray*}
		n &\geq& \frac{(b - a)^2}{2 \varepsilon^2} \log\left(\frac{2}{\delta}\right),
	\end{eqnarray*}
	then it holds that
	\begin{eqnarray*}
		\pr\left(\left|\bar{x}-\E(X)\right|\geq\epsilon\right) &\leq&  \delta.
	\end{eqnarray*}
\end{corollary}

\subsection{Notations}
For a matrix $G \in \reals^{m \times n}$, we use $G^{\top}$, $G^{\dagger}$, $\range(G)$, $\zspace(G)$, $\|G\|_2$, and $\|G\|_F$ to denote the transpose, the Moore-Penrose pseudoinverse, the column space, the null space, the spectral norm, and the Frobenius norm of $G$. For two Hermitian matrices $G$ and $H$ of the same size, $G\succ H$ means that $G-H$ is positive definite, where $\succ$ denotes the strict Loewner order. If the matrix $G$ is symmetric and positive semi-definite, then $\lambda_{\min}(G)$ and  $\lambda_{\max}(G)$ represent the smallest and the largest positive eigenvalues of $G$. For a general matrix $G$, we denote by $\sigma_{i}(G)$ its positive singular values, and in particular, $\sigma_{\min}(G)$ and $\sigma_{\max}(G)$ the smallest and largest positive singular values. We denote by $\mathbb{DR}^{m\times m}$ the set of all invertible diagonal matrices in $\reals^{m\times m}$.
Moreover, $x_{\star}$ denotes the unique least-norm solution of (\ref{linear-system}), $r_j^{(k)}=b_j-A^{(j)}x^{(k)}$ the residual of the $j$th equation of (\ref{linear-system}) at $x^{(k)}$, $[m]$ the index set including integers $\{1,\ldots,m\}$, $|\cdot|$ the absolute value of a real number or the number of elements in a set.

In addition, we denote by $\mathbb{E}_k$ the conditional expectation given the first $k$ iterations, that is,
\begin{eqnarray*}
	\mathbb{E}_k(\cdot)=\mathbb{E}\left(\cdot \mid \tau^{(0)}, \tau^{(1)}, \ldots, \tau^{(k-1)}\right),
\end{eqnarray*}
where $\tau^{(\ell)}$ is denoted by the index of the rows selected at the $\ell$th iterate for $\ell=0,1, \ldots, k-1$. According to the law of iterated expectation, we have  $\mathbb{E}\left(\mathbb{E}_k(\cdot)\right)=\mathbb{E}(\cdot)$.

\section{The RBSK Method}
\label{sec:RBSK}

\begin{definition}\label{def:batchSampling}
  A random vector $\tau$ is called a \textsc{batch-sampling} of the index set $[m]$ of batch-size $q$, if
  \begin{enumerate}
    \item $\tau=(\tau_1,\tau_2,\ldots,\tau_q)$ with $\tau_i\in [m]$ satisfies the joint distribution $\bbb{P}$, i.e.,
  \begin{eqnarray}
  \label{jointDistrib}
    \pr(\tau=(j_1,\ldots,j_q)) &=& \bbb{p}_{j_1,\ldots,j_q}\geq 0\ 
    \mbox{ with }\sum_{j_1,\ldots,j_q=1}^{m}\bbb{p}_{j_1,\ldots,j_q}=1,
  \end{eqnarray}
    \item and $\tau_i$ are random variables satisfying the following marginal distributions
  \begin{eqnarray}
  \label{marDistrib} \nonumber 
    \pr(\tau_i=j) &=& \sum_{j_1,\ldots,j_{i-1},j_{i+1},\ldots,j_q=1}^{m}\bbb{p}_{j_1,\ldots,j_{i-1},j_i=j,j_{i+1},\ldots,j_q} \\
    &\triangleq& p_{ij}\geq 0\ 
    \mbox{ with } \sum_{j=1}^{m} p_{ij}=1, \mbox{ for } i=1,2,\dots,q,
  \end{eqnarray}
    \item and the corresponding diagonal matrices
  \begin{eqnarray}
  \label{marDistribMatrix}
  P_i=\diag\{p_{i1},p_{i2},\ldots,p_{im}\},
        \ \mbox{ for } i=1,\ldots,q,
  \end{eqnarray}
  admit the property below
  \begin{eqnarray}
  \label{marDistribMatrixProperty}
    \sum_{i=1}^{q} P_i &\succ& 0.
  \end{eqnarray}
  \end{enumerate}
\end{definition}

The joint distribution $\bbb{P}$ defined by (\ref{jointDistrib}) can be naturally represented by a $q$th-order tensor, which plays the role of a parameter in defining an instance of randomized block Kaczmarz methods. A different $\bbb{P}$ leads to a different stochastic sampling (i.e., drawing batch-sampling $\tau$ in a different way), which in turn defines a different randomized block Kaczmarz instance.

As a matter of fact, the indices in $[m]$ may be even allowed to appear repeatedly in one draw of the batch-sampling $\tau$ for some prescribed joint distributions $\bbb{P}$. If we denote by $n_{\txttiny{E}}(\tau)$ the number of the effective row indices included in $\tau$ (i.e., the number of unique row indices in $\tau$), it satisfies that
\begin{eqnarray*}
  1 \leq n_{\txttiny{E}}(\tau) \leq q.
\end{eqnarray*}
Obviously, when $n_{\txttiny{E}}(\tau)<q$, there are repeated row indices appearing in $\tau$. Therefore, we denote by $\tau_{\txttiny{E}}$ the effective set of the batch-sampling $\tau$ (i.e., the set of the effective row indices included in $\tau$). Then, the range of $\tau$ can be defined by the set of all possible values of the effective set of the batch-sampling $\tau$ (i.e., the set of all possible values of $\tau_{\txttiny{E}}$), which is denoted by $\range(\tau)$. Since different draws of the batch-sampling $\tau$ may lead to a same effective set $\tau_{\txttiny{E}}$, for a prescribed effective set $\tau_{\txttiny{E}}$, we define the preimage of $\tau_{\txttiny{E}}$ as below
\begin{eqnarray*}
  \tau^{-1}\left(\tau_{\txttiny{E}}\right) &=& \left\{
  \widehat{\tau}\ |\ \mbox{the effective set of $\widehat{\tau}$ is $\tau_{\txttiny{E}}$}
  \right\},
\end{eqnarray*}
i.e., the set of all possible draws of the batch-sampling $\tau$ that lead to the same effective set $\tau_{\txttiny{E}}$. Thanks to the above notations, we also define the probability that one draw of the batch-sampling $\tau$ leads to a prescribed effective set $\tau_{\txttiny{E}}$ as
\begin{eqnarray*}
  \pr\left(\tau_{\txttiny{E}}\right) &=& \pr\left(\tau\in\tau^{-1}\left(\tau_{\txttiny{E}}\right)\right) \\
  &=& \sum_{\tau\in\tau^{-1}\left(\tau_{\txttiny{E}}\right)} \pr(\tau=(j_1,\ldots,j_q)).
\end{eqnarray*}
For instance, let the data matrix $A\in\reals^{4\times n}$, the batch-size $q=3$
(i.e., the batch-sampling $\tau=(\tau_1,\tau_2,\tau_3)$), and a prescribed effective set $\tau_{\txttiny{E}}=\{1,2\}$, then the preimage of $\tau_{\txttiny{E}}$ reads
\begin{eqnarray*}
  \tau^{-1}\left(\tau_{\txttiny{E}}\right) &=& \left\{\tau=(1,1,2),\tau=(1,2,1),\tau=(2,1,1),\right. \\
  & & \left.\tau=(2,2,1),\tau=(2,1,2),\tau=(1,2,2)\right\},
\end{eqnarray*}
and the probability that one draw of the batch-sampling $\tau=(\tau_1,\tau_2,\tau_3)$ leads to the effective set $\tau_{\txttiny{E}}=\{1,2\}$ is given below
\begin{eqnarray*}
  \pr\left(\tau_{\txttiny{E}}=\{1,2\}\right) &=& \pr(\tau=(1,1,2))+\pr(\tau=(1,2,1))+\pr(\tau=(2,1,1)) \\
  & & +\pr(\tau=(2,2,1))+\pr(\tau=(2,1,2))+\pr(\tau=(1,2,2)).
\end{eqnarray*}

In order to describe the RBSK method, we denote by $G_{\tau}$ the batch-sampling matrix
including the rows of a matrix $G\in\reals^{m\times n}$ indexed by $\tau$, $x_{\tau}$ the batch-sampling vector including the entries of a vector $x\in\reals^m$ indexed by $\tau$. The RBSK method (Algorithm \ref{alg-RBSK}) is interpreted as follows.

\begin{breakalgo}{The RBSK method}{alg-RBSK}
\begin{algorithmic}[1]
  \State \textbf{Input}: Initial guess $x^{(0)}$, stopping criterion, joint distribution $\bbb{P}$, maximal iteration count MaxIter;
  \State \textbf{Output}: approximate solution $x$;
  \For{$k = 0$ to MaxIter-1}
    \If{the stopping criterion holds}
        \State $x=x^{(k)}$; \textbf{break};
    \EndIf
    \State Draw a batch-sampling $\tau^{(k)}\sim\bbb{P}$;
    \State Project the $k$th iterate $x^{(k)}$ onto the solution space of $A_{\tau^{(k)}}x=b_{\tau^{(k)}}$ by (\ref{BRK-method});
  \EndFor
\end{algorithmic}
\end{breakalgo}

In Algorithm \ref{alg-RBSK}, the sampling rules controlled by the batch-sampling $\tau\sim\bbb{P}$ may recover some existing sampling rules (see the three examples in Appendix \ref{examples_samplingRules}), and also lead to new ones.

\section{Main Lemmas}
\label{sec:lemma}

To study the convergence property of the RBSK method, the following main lemmas are needed. All other relevant lemmas are presented in Appendix \ref{otherLemmas}.

\begin{lemma}
  \label{lem6}
  Let $\tau\sim\bbb{P}$ be a batch-sampling from the index set $[m]$ of
  batch-size $q$, whose entries $\tau_i$ are random
  variables satisfying the marginal distributions {\rm(\ref{marDistrib})},
  $P_i$ be the diagonal matrices defined in {\rm(\ref{marDistribMatrix})}. Let $A\in\reals^{m\times n}$ be a
  nonzero matrix, $r\in\range(A)$ be a nonzero vector, and $S\in\mathbb{DR}^{m\times m}$ be a nonsingular diagonal matrix. Let $A_{\tau}$ be the batch-sampling matrix including the rows of $A$ indexed by $\tau$, $r_{\tau}$ be the batch-sampling vector including the entries of $r$ indexed by $\tau$, and $S_{\tau}$ be the batch-sampling principal submatrix including the diagonals of $S$ indexed by $\tau$. Then it holds that
  \begin{eqnarray}
  \label{lem5Eq0}
    \E\left(\|A_{\tau}^{\dag}r_{\tau}\|_2^2\right) &\geq& \sum_{i=1}^{q}\|B_{S;i}^{-\frac{1}{2}} P_i^{\frac{1}{2}}S r\|_2^2,
  \end{eqnarray}
  where
  \begin{eqnarray}
  \label{lem5Eq3}
    B_{S;i}=\begin{bmatrix}
     \beta^S_{i,1} &             &        &           \\
                 & \beta^S_{i,2} &        &           \\
                 &             & \ddots &           \\
                 &             &        & \beta^S_{i,m}
        \end{bmatrix},\quad\mbox{for $i=1,\ldots,q$,}
  \end{eqnarray}
  with
  \begin{eqnarray*}
    \beta^S_{i,j} &=& \left\{\begin{array}{lr}
    \max_{\tau=(\tau_1,\ldots,\tau_{i-1},\tau_i=j,
    \tau_{i+1},\ldots,\tau_q)} \|S_{\tau} A_{\tau}\|_2^2, & \mbox{if $p_{ij}>0$}, \\
        \max_{\tau\sim\bbb{P}} \|S_{\tau} A_{\tau}\|_2^2, & \mbox{if $p_{ij}=0$},
                             \end{array}
    \right.
    \quad\mbox{for $j=1,\ldots,m$}.
  \end{eqnarray*}
\end{lemma}
{\em Proof.}
Since $r\in\range(A)$, it holds that $r_{\tau}\in\range(A_{\tau})$. In addition, since $S$ is nonsingular, so is $S_{\tau}$.
According to Lemmas \ref{lem2} and \ref{lem3}, it follows
that
\begin{eqnarray*}
  \|A_{\tau}^{\dag}r_{\tau}\|_2^2 &=& \|(S_{\tau}A_{\tau})^{\dag}S_{\tau}r_{\tau}\|_2^2 \\
  &\geq& \frac{1}{\|S_{\tau}A_{\tau}\|_2^2}\|S_{\tau}r_{\tau}\|_2^2.
\end{eqnarray*}
Due to the order preserving of expectation, it reads that
\begin{eqnarray}
 \nonumber 
  \E\left(\|A_{\tau}^{\dag}r_{\tau}\|_2^2\right) &\geq& \E\left(\frac{1}{\|S_{\tau}A_{\tau}\|_2^2}\|S_{\tau}r_{\tau}\|_2^2\right) \\ \nonumber
  &=& \E\left(\sum_{i=1}^{q}|s_{\tau_i}|^2|r_{\tau_i}|^2\frac{1}{\|S_{\tau}A_{\tau}\|_2^2}\right) \\ \label{lem5Eq1}
  &=& \sum_{i=1}^{q}\E\left(|s_{\tau_i}|^2|r_{\tau_i}|^2\frac{1}{\|S_{\tau}A_{\tau}\|_2^2}\right).
\end{eqnarray}
According to Lemma \ref{lem4},
for each $\E\left(|s_{\tau_i}|^2|r_{\tau_i}|^2\frac{1}{\|S_{\tau}A_{\tau}\|_2^2}\right)$ in (\ref{lem5Eq1}), it holds that
\begin{eqnarray}
 \nonumber 
  \E\left(|s_{\tau_i}|^2|r_{\tau_i}|^2\frac{1}{\|S_{\tau}A_{\tau}\|_2^2}\right) &=& \E\left[|s_{\tau_i}|^2|r_{\tau_i}|^2 \E\left(\left.\frac{1}{\|S_{\tau}A_{\tau}\|_2^2}\right|_{|s_{\tau_i}|^2|r_{\tau_i}|^2}\right)\right] \\ \nonumber
  &=& \E\left[|s_{\tau_i}|^2|r_{\tau_i}|^2 \E\left(\left.\frac{1}{\|S_{\tau}A_{\tau}\|_2^2}\right|_{\tau_i}\right)\right] \\ \nonumber
  &=& \sum_{j=1}^{m} |s_j|^2|r_j|^2 p_{ij} \ \E\left(\left.\frac{1}{\|S_{\tau}A_{\tau}\|_2^2}\right|_{\tau_i=j}\right) \\ \nonumber
  &\geq& \sum_{j=1}^{m} |s_j|^2|r_j|^2 p_{ij} \frac{1}{\beta^S_{i,j}} \\ \label{lem5Eq2}
  &=& \|B_{S;i}^{-\frac{1}{2}} P_i^{\frac{1}{2}}S r\|_2^2,
\end{eqnarray}
where $p_{ij}\geq 0$ are marginal distributions for random variables $\tau_i$ introduced in Definition \ref{def:batchSampling}.
Obviously, (\ref{lem5Eq1}) and (\ref{lem5Eq2}) lead to (\ref{lem5Eq0}). $\hfill\square$

\begin{lemma}
  \label{lem7}
  Let $A\in\reals^{m\times n}$ be a matrix with $\rank(A)=d>0$. Let $\tau\sim\bbb{P}$ be a batch-sampling
  from the index set $[m]$ of batch-size $q$, whose entries $\tau_i$ are random
  variables satisfying the marginal distributions {\rm(\ref{marDistrib})}. Let $A_{\tau}$ be a batch-sampling matrix including the rows of $A$ indexed by $\tau$, and $d_{\tau}=\dim\left[\range(A_{\tau}^{\T})\right]$. Let $\mathcal{W}_{\tau}\subseteq\range(A^{\T})$ be a subspace satisfying the facts
  $\mathcal{W}_{\tau}\perp\range(A_{\tau}^{\T})$ and $\dim(\mathcal{W}_{\tau})=d-d_{\tau}$.
  Let the columns of matrices  $V_{\tau}\in\reals^{n\times d_{\tau}}$ and $W_{\tau}\in\reals^{n\times (d-d_{\tau})}$
  be the orthonormal bases of the subspaces $\range(A_{\tau}^{\T})$ and $\mathcal{W}_{\tau}$, respectively. If the vector $y\in\range(A^{\T})$ satisfies
  $y\perp \range(A_{\tau}^{\T})$, then
  \begin{eqnarray*}
    \|DAy\|_2^2 &\geq& \xi_{\tau}\|y\|_2^2 \quad \mbox{with $\xi_{\tau}=\lambda_{\min}\left( W_{\tau}^{\T}A^{\T}D^{\T}DAW_{\tau} \right)$},
  \end{eqnarray*}
  for all  invertible matrices $D\in\reals^{m\times m}$.
\end{lemma}
{\em Proof.}
Due to the nonsingularity of $D$, and the definitions of $V_{\tau}$ and $W_{\tau}$, it reads that
\begin{eqnarray*}
  \range\left(A^{\T}\right) = \range\left[(DA)^{\T}\right] = \range\left\{[V_{\tau}\ W_{\tau}]\right\} \quad \mbox{and $V_{\tau}^{\T} W_{\tau}=0$}.
\end{eqnarray*}
Thanks to the above facts, there exists a full column rank matrix $C\in\reals^{m\times d}$ such that
\begin{eqnarray*}
  (DA)^{\T} &=& [V_{\tau}\ W_{\tau}]\ C^{\T},
\end{eqnarray*}
thus it follows that
\begin{eqnarray*}
  DAW_{\tau} &=& C\ [V_{\tau}\
                      W_{\tau}]^{\T} W_{\tau} \\
             &=& C\ [0\ I]^{\T} \\
             &=& C_{-d_{\tau}},
\end{eqnarray*}
where $C_{-d_{\tau}}\in\reals^{m\times (d-d_{\tau})}$ represents the submatrix obtained by deleting the first $d_{\tau}$ columns of $C$. Since $C$ is
a full column rank matrix, so is $C_{-d_{\tau}}$, which leads to a fact that
\begin{eqnarray*}
  (DAW_{\tau})^{\T}DAW_{\tau} &=& W_{\tau}^{\T}A^{\T}D^{\T}DAW_{\tau} \\
  &=& C_{-d_{\tau}}^{\T}C_{-d_{\tau}}
\end{eqnarray*}
is a real symmetric positive definite matrix.

If a vector $u\in\mathcal{W}_{\tau}$ satisfies $\|u\|_2=1$, then there exists a vector $z\in\reals^{d-d_{\tau}}$
such that $u=W_{\tau}z$ and $\|z\|_2=1$. Thus, a constant $\xi_{\tau}$ may be defined as
\begin{eqnarray*}
  \xi_{\tau} &=& \min_{u\in\mathcal{W}_{\tau},\ \|u\|_2=1} \|DAu\|_2^2 \\
  &=& \min_{z\in\reals^{d-d_{\tau}},\ \|z\|_2=1} \|DAW_{\tau}z\|_2^2 \\
  &=& \lambda_{\min} \left(W_{\tau}^{\T}A^{\T}D^{\T}DAW_{\tau}\right).
\end{eqnarray*}
In addition, since the vector $y\in\range(A^{\T})$ satisfies
$y\perp \range(A_{\tau}^{\T})$, it holds that $y\in\mathcal{W}_{\tau}$, which leads to
\begin{eqnarray*}
  \|DAy\|_2^2 &=& \left\|DA\frac{y}{\|y\|_2}\right\|_2^2 \|y\|_2^2 \\
  &\geq& \left(\min_{u\in\mathcal{W}_{\tau},\ \|u\|_2=1} \|DAu\|_2^2\right) \|y\|_2^2 \\
  &=& \xi_{\tau}\|y\|_2^2.
\end{eqnarray*}
$\hfill\square$

\begin{lemma}
  \label{lem8}
  Let $\tau\sim\bbb{P}$ be a batch-sampling from the index set $[m]$ of batch-size $q$.
  If $j\in [m]$ is a row index, and $\tau_{\txttiny{E$(j)$}}\in\range(\tau)$ are effective sets
  that include the row index $j$ (i.e., $j\in\tau_{\txttiny{E$(j)$}}$), then
  \begin{eqnarray}
  \label{lem7EqResult}
    \sum_{\tau_{\txttiny{E$(j)$}}\in\range(\tau)} \pr\left(\tau_{\txttiny{E$(j)$}}\right) \leq
    \sum_{i=1}^{q} p_{ij} \leq
    \sum_{\tau_{\txttiny{E$(j)$}}\in\range(\tau)} \left(q-|\tau_{\txttiny{E$(j)$}}|+1\right)\pr\left(\tau_{\txttiny{E$(j)$}}\right),\
    \mbox{for }\ j=1,\ldots,m,
  \end{eqnarray}
  where $p_{ij}=\pr(\tau_i=j)$ ($i=1,\ldots,q$) are the marginal distributions of random variables $\tau_i$
  defined by {\rm(\ref{marDistrib})}, and $|\tau_{\txttiny{E$(j)$}}|$ is the number of unique row indices in the effective set $\tau_{\txttiny{E$(j)$}}$. In addition, the equalities in {\rm(\ref{lem7EqResult})} get satisfied if $|\tau_{\txttiny{E$(j)$}}|=q$.
\end{lemma}
{\em Proof.}
Let $\widehat{\tau}$ be one draw of the batch-sampling $\tau$ such that
\begin{eqnarray*}
  \widehat{\tau} &\in& \tau^{-1}\left(\tau_{\txttiny{E$(j)$}}\right),
\end{eqnarray*}
then it holds that
\begin{eqnarray}
 \nonumber 
  \sum_{i=1}^{q} p_{ij} &=& \sum_{i=1}^{q}\
  \sum_{j_1,\ldots,j_{i-1},j_{i+1},\ldots,j_q=1}^{m}\bbb{p}_{j_1,\ldots,j_{i-1},j_i=j,j_{i+1},\ldots,j_q} \\ \label{lem7Eq1}
  &=& \sum_{\tau_{\txttiny{E$(j)$}}\in\range(\tau)}\ \sum_{\widehat{\tau}\in\tau^{-1}\left(\tau_{\txttiny{E$(j)$}}\right)}
  n\left(\widehat{\tau},j\right)\pr(\tau=\widehat{\tau}),
\end{eqnarray}
where $n\left(\widehat{\tau},j\right)$ represents the number of repeats of the row index $j$ in $\widehat{\tau}$,
and $n\left(\widehat{\tau},j\right)$ satisfies
\begin{eqnarray}
\label{lem7Eq2}
  1 \leq n\left(\widehat{\tau},j\right) \leq \left(q-|\tau_{\txttiny{E$(j)$}}|+1\right).
\end{eqnarray}
The equality (\ref{lem7Eq1}) is due to two facts: firstly, if the row index $j$ repeats $n\left(\widehat{\tau},j\right)$
times in $\widehat{\tau}$, the probability $\pr(\tau=\widehat{\tau})$ repeats $n\left(\widehat{\tau},j\right)$ times
in $\sum_{i=1}^{q} p_{ij}$; secondly, if the effective sets $\tau_{\txttiny{E}}$ and
$\widetilde{\tau}_{\txttiny{E}}$ are not equal
(i.e., $\tau_{\txttiny{E}}\neq\widetilde{\tau}_{\txttiny{E}}$), the intersection of the corresponding preimages
is empty (i.e.,
$\tau^{-1}\left(\tau_{\txttiny{E}}\right)\cap\tau^{-1}\left(\widetilde{\tau}_{\txttiny{E}}\right)=\emptyset$).

Obviously, the equality (\ref{lem7Eq1}) and the inequality (\ref{lem7Eq2}) lead to
\begin{eqnarray*}
  \sum_{\tau_{\txttiny{E$(j)$}}\in\range(\tau)}\ \sum_{\widehat{\tau}\in\tau^{-1}\left(\tau_{\txttiny{E$(j)$}}\right)}
  \pr(\tau=\widehat{\tau})
  \leq \sum_{i=1}^{q} p_{ij} \leq
  \sum_{\tau_{\txttiny{E$(j)$}}\in\range(\tau)}\left(q-|\tau_{\txttiny{E$(j)$}}|+1\right)
  \sum_{\widehat{\tau}\in\tau^{-1}\left(\tau_{\txttiny{E$(j)$}}\right)}
  \pr(\tau=\widehat{\tau}).
\end{eqnarray*}
The above inequality together with the fact
\begin{eqnarray*}
  \pr\left(\tau_{\txttiny{E$(j)$}}\right) &=& \sum_{\widehat{\tau}\in\tau^{-1}\left(\tau_{\txttiny{E$(j)$}}\right)}
  \pr(\tau=\widehat{\tau})
\end{eqnarray*}
leads to the inequality (\ref{lem7EqResult}).
$\hfill\square$

\section{Main Theorems}
\label{sec:convAnal}
In this section, we present the convergence analysis of RBSK. The main result is Theorem~\ref{convRBSK}, which provides a new convergence rate bound for the RBSK method. Remarks~\ref{remark5.1}--\ref{remark5.3} clarify the respective roles of the scaling matrix $S$ and the joint distribution $\bbb{P}$ in sharpening the bound and determining the iteration of RBSK, which also explain the scale-invariant property of the new bound. Corollary~\ref{convRBSK_relaxed} then presents a relaxed version of Theorem~\ref{convRBSK} in a more concise form. Remark~\ref{remark:relaxOne} specifies this relaxed bound for several representative batch-sampling rules, and proves that even the relaxed bound is sharper than the classical bound in \cite{NeedellTropp2014LAA}. Motivated by Corollary~\ref{convRBSK_relaxed}, Theorem~\ref{convRBSK_relaxed_pave} gives an improved bound for the row-paving case, while Remark~\ref{remark:relaxTwo} further analyzes the improvement of Theorem~\ref{convRBSK_relaxed_pave} over Corollary~\ref{convRBSK_relaxed} and the classical bound in \cite{NeedellTropp2014LAA}. Finally, under slightly stronger assumptions, Theorem~\ref{convRBSK_sharper} derives an even sharper estimate than Theorem~\ref{convRBSK} by using concentration inequalities. We now state our main result.

\begin{theorem}
  \label{convRBSK}
  Let the linear system {\rm(\ref{linear-system})} be consistent. Let $A\in\reals^{m\times n}$ be a matrix with $\rank(A)=d>0$. Let $\tau\sim\bbb{P}$ be a batch-sampling from the index set $[m]$ of batch-size $q$, whose entries $\tau_i$ are random
  variables satisfying the marginal distributions {\rm(\ref{marDistrib})}. Let $A_{\tau}$ be a batch-sampling matrix including the rows of $A$ indexed by $\tau$, and $d_{\tau}=\dim\left[\range(A_{\tau}^{\T})\right]$. Let $\mathcal{W}_{\tau}\subseteq\range(A^{\T})$ be a subspace satisfying the facts
  $\mathcal{W}_{\tau}\perp\range(A_{\tau}^{\T})$ and $\dim(\mathcal{W}_{\tau})=d-d_{\tau}$.
  Let the columns of matrix $W_{\tau}\in\reals^{n\times (d-d_{\tau})}$ be the orthonormal basis of the subspace $\mathcal{W}_{\tau}$.
  If the initial guess $x^{(0)}$
  in the RBSK method satisfies $x^{(0)}\in\range(A^{\T})$,
  then the iterative sequence
  $\{x^{(k)}\}_{k=0}^{+\infty}$ converges to the least-norm solution $x_{\star}$ of {\rm(\ref{linear-system})} in
  expectation. In addition, the mean squared error satisfies
  \begin{eqnarray}
    \label{thm1Eq7}
    \E\left(\|x^{(k+1)}-x_{\star}\|_2^2\right)
    &\leq& \left[\min_{S\in\mathbb{DR}^{m\times m}}\left(1-\xi\right)\right]^k
    \left[\min_{S\in\mathbb{DR}^{m\times m}}\left(1-\eta\right)\right]\|x^{(0)}-x_{\star}\|_2^2,
  \end{eqnarray}
  in particular, for $k=0$, it satisfies
  \begin{eqnarray}
    \label{thm1Eq2}
    \E\left(\|x^{(1)}-x_{\star}\|_2^2\right) 
    &\leq& \left[\min_{S\in\mathbb{DR}^{m\times m}}\left(1-\eta\right)\right]\|x^{(0)}-x_{\star}\|_2^2,
  \end{eqnarray}
  and, for $k>0$, it satisfies
  \begin{eqnarray}
    \label{thm1Eq5}
      \E\left(\|x^{(k+1)}-x_{\star}\|_2^2\right)
      &\leq& \left[\min_{S\in\mathbb{DR}^{m\times m}}\left(1-\xi\right)\right]\E\left(\|x^{(k)}-x_{\star}\|_2^2\right),
  \end{eqnarray}
  where
  \begin{eqnarray*}
    \eta &=&\lambda_{\min}\left(A^{\T}D^2A\right),
  \end{eqnarray*}
  and
  \begin{eqnarray*}
    \xi &=& \min_{\tau\sim\bbb{P}} \lambda_{\min}\left(W_{\tau}^{\T}A^{\T} D^2AW_{\tau}\right)
  \end{eqnarray*}
  with
  \begin{eqnarray*}
  D^2 &=& S^{\T}\left(\sum_{i=1}^{q} B_{S;i}^{-1} P_i\right)S,
\end{eqnarray*}
  here, $S\in\mathbb{DR}^{m\times m}$ is any prescribed nonsingular diagonal matrix, $P_i\in\mathbb{R}^{m\times m}$ and $B_{S;i}\in\mathbb{R}^{m\times m}$ are diagonal matrices given in {\rm(\ref{marDistribMatrix})} and {\rm(\ref{lem5Eq3})}.
\end{theorem}
{\em Proof.}
As stated in \cite{BaiWu18SIAM}, if the RBSK method converges to a solution to the consistent linear system (\ref{linear-system}), and the initial guess $x^{(0)}$ of the iteration belongs to the column space of $A^{\T}$, then this solution must be the least-norm solution $x_{\star}=A^{\dag}b$.

According to the definition of the RBSK method, the $k$th iterate $x^{(k)}$ is the orthogonal projection of $x^{(k-1)}$ onto the solution space of the $(k-1)$th sampled linear system $A_{\tau^{(k-1)}}x=b_{\tau^{(k-1)}}$, which leads to a fact that
\begin{eqnarray*}
  r_{\tau^{(k-1)}}^{(k)} = b_{\tau^{(k-1)}}-A_{\tau^{(k-1)}}x^{(k)}=0.
\end{eqnarray*}
Similarly to $x^{(k)}$, the next iterate $x^{(k+1)}$ is the orthogonal projection of $x^{(k)}$ onto the solution space of the $k$th sampled linear system $A_{\tau^{(k)}}x=b_{\tau^{(k)}}$, thus the vectors $x^{(k+1)}-x_{\star}$, $x^{(k)}-x_{\star}$, and $x^{(k+1)}-x^{(k)}$ satisfy the following fact
\begin{eqnarray*}
  \|x^{(k)}-x_{\star}\|_2^2 &=& \|x^{(k+1)}-x_{\star}\|_2^2 +\|x^{(k+1)}-x^{(k)}\|_2^2,
\end{eqnarray*}
equivalently, it reads that
\begin{eqnarray*}
  \|x^{(k+1)}-x_{\star}\|_2^2 &=& \|x^{(k)}-x_{\star}\|_2^2 - \|x^{(k+1)}-x^{(k)}\|_2^2  \\
  &=& \|x^{(k)}-x_{\star}\|_2^2 - \|A_{\tau^{(k)}}^{\dag}r_{\tau^{(k)}}^{(k)}\|_2^2.
\end{eqnarray*}

Taking conditional expectation conditioned on the previous $k$ iterations of the RBSK method, making use of the fact $r^{(k)}=b-Ax^{(k)}\in\range(A)$, and applying Lemma \ref{lem6}, it holds that
\begin{eqnarray}
 \nonumber 
  \E_k\left(\|x^{(k+1)}-x_{\star}\|_2^2\right) &=& \|x^{(k)}-x_{\star}\|_2^2 - \E_k\left(\|A_{\tau^{(k)}}^{\dag}r_{\tau^{(k)}}^{(k)}\|_2^2\right) \\ \nonumber
  &\leq& \|x^{(k)}-x_{\star}\|_2^2 - \sum_{i=1}^{q}\|B_{S;i}^{-\frac{1}{2}} P_i^{\frac{1}{2}}Sr^{(k)}\|_2^2 \\ \nonumber
  &=& \|x^{(k)}-x_{\star}\|_2^2 - \sum_{i=1}^{q}\|B_{S;i}^{-\frac{1}{2}} P_i^{\frac{1}{2}}SA(x^{(k)}-x_{\star})\|_2^2 \\ \nonumber
  &=& \|x^{(k)}-x_{\star}\|_2^2 - \sum_{i=1}^{q} (x^{(k)}-x_{\star})^{\T}A^{\T}(B_{S;i}^{-\frac{1}{2}} P_i^{\frac{1}{2}}S)^2A(x^{(k)}-x_{\star}) \\ \nonumber
  &=& \|x^{(k)}-x_{\star}\|_2^2 - (x^{(k)}-x_{\star})^{\T}A^{\T}S^{\T}\left(\sum_{i=1}^{q}B_{S;i}^{-1} P_i\right)SA(x^{(k)}-x_{\star}) \\ \label{thm1Eq1}
  &=& \|x^{(k)}-x_{\star}\|_2^2 - \|DA(x^{(k)}-x_{\star})\|_2^2,
\end{eqnarray}
where the diagonal matrix
\begin{eqnarray}
\label{D_formula}
  D &=& \left[S^{\T}\left(\sum_{i=1}^{q}B_{S;i}^{-1} P_i\right)S\right]^{\frac{1}{2}}\succ 0
\end{eqnarray}
is well defined, since $S$ is a nonsingular diagonal real matrix, $B_{S;i}$ is diagonally positive definite,
and the matrices $P_i$ are diagonally positive semi-definite and satisfy the property (\ref{marDistribMatrixProperty}).

Since the initial guess $x^{(0)}\in\range(A^{\T})$ and the correction $A_{\tau^{(k)}}^{\dag}r_{\tau^{(k)}}^{(k)}$ adopted to update $x^{(k)}$ at each iteration
guarantee $x^{(k)}\in\range(A^{\T})$ for all $k\geq 0$, together with the fact that the least-norm solution satisfies $x_{\star}\in\range(A^{\T})$,
it can be concluded that $x^{(k)}-x_{\star}\in\range(A^{\T})$.

For $k=0$, since $x^{(1)}-x_{\star}\in\range(A^{\T})$ leads to a fact $x^{(1)}-x_{\star}\in\range[(DA)^{\T}]$ due to the nonsingularity of the diagonal matrix $D$, the relation (\ref{thm1Eq1}) together with Lemma \ref{lem1} leads to a
convergence rate estimate at the 1st iteration as follows
\begin{eqnarray*}
  \E\left(\|x^{(1)}-x_{\star}\|_2^2\right) 
  &\leq& \left(1-\eta\right)\|x^{(0)}-x_{\star}\|_2^2
\end{eqnarray*}
with
\begin{eqnarray*}
  \eta&=&\lambda_{\min}\left(A^{\T}D^2A\right),
\end{eqnarray*}
thanks to the fact that $S\in\mathbb{DR}^{m\times m}$ is an arbitrary nonsingular diagonal matrix,
the above estimate can be minimized with respect to $S$, which results in the estimate (\ref{thm1Eq2}).

For $k>0$, due to the fact $r_{\tau^{(k-1)}}^{(k)}=0$, i.e., $A_{\tau^{(k-1)}}(x^{(k)}-x_{\star})=0$,
it reads that $x^{(k)}-x_{\star} \perp \range(A_{\tau^{(k-1)}}^{\T})$.
According to Lemma \ref{lem7}, since $x^{(k)}-x_{\star}\in\range(A^{\T})$ and $x^{(k)}-x_{\star} \perp \range(A_{\tau^{(k-1)}}^{\T})$, the relation (\ref{thm1Eq1}) leads to the following convergence rate estimate at the $k$th iteration
\begin{eqnarray}
 \label{thm1Eq3}
  \E_k\left(\|x^{(k+1)}-x_{\star}\|_2^2\right) 
  &\leq& \left(1-\xi_{\tau^{(k-1)}}\right)\|x^{(k)}-x_{\star}\|_2^2
\end{eqnarray}
with
\begin{eqnarray*}
  \xi_{\tau^{(k-1)}} &=& \lambda_{\min}\left(W_{\tau^{(k-1)}}^{\T}A^{\T} D^2AW_{\tau^{(k-1)}}\right).
\end{eqnarray*}
By defining a constant $\xi$ independent of $k$ as follows
\begin{eqnarray*}
  \xi &=& \min_{\tau\sim\bbb{P}} \lambda_{\min}\left(W_{\tau}^{\T}A^{\T} D^2AW_{\tau}\right),
\end{eqnarray*}
the estimate (\ref{thm1Eq3}) leads to
\begin{eqnarray}
\label{thm1Eq4}
  \E_k\left(\|x^{(k+1)}-x_{\star}\|_2^2\right)
  &\leq& \left(1-\xi\right)\|x^{(k)}-x_{\star}\|_2^2.
\end{eqnarray}
By taking full expectation on both sides of (\ref{thm1Eq4}), it holds that
\begin{eqnarray*}
  \E\left(\|x^{(k+1)}-x_{\star}\|_2^2\right)
  &\leq& \left(1-\xi\right)\E\left(\|x^{(k)}-x_{\star}\|_2^2\right),
\end{eqnarray*}
due to the arbitrariness of $S$, the above inequality leads to the estimate (\ref{thm1Eq5}).
Combining (\ref{thm1Eq2}) and (\ref{thm1Eq5}), and through recursion with respect to $k$, the estimate (\ref{thm1Eq7})
is obtained.
$\hfill\square$

\begin{remark}
	\label{remark5.1}
  According to {\rm(\ref{thm1Eq5})}, for a prescribed joint distribution $\bbb{P}$ and iteration count $k>0$,
  the convergence rate of the mean squared error is bounded by
  \begin{eqnarray}
    U(\bbb{P},S) &=& 1-\min_{\tau\sim\bbb{P}} \lambda_{\min}\left(W_{\tau}^{\T}A^{\T} D^2AW_{\tau}\right).
  \end{eqnarray}
  The nonsingular diagonal matrix does not serve as an iteration parameter for the RBSK method. In fact,
  due to the arbitrariness of $S$, it serves as a parameter to minimize the upper bound $U(\bbb{P},S)$,
  or equivalently, to achieve the sharpest upper bound of the convergence rate of the RBSK method
  when the joint distribution $\bbb{P}$ is given, i.e.,
  \begin{eqnarray*}
    \widehat{S}_{\txt{opt}} &=& \arg\min_{S\in\mathbb{DR}^{m\times m}} U(\bbb{P},S),
  \end{eqnarray*}
  which leads to
  \begin{eqnarray*}
    U(\bbb{P},\widehat{S}_{\txt{opt}}) &=& \min_{S\in\mathbb{DR}^{m\times m}} U(\bbb{P},S).
  \end{eqnarray*}
  Therefore, a specific value of $S$ only affects the value of $U(\bbb{P},S)$ rather than the actual
  convergence behavior of the RBSK method.
\end{remark}

\begin{remark}
\label{remark5.2}
  Different from the role of the parameter $S$, the joint distribution $\bbb{P}$ serves as
  an iteration parameter for the RBSK method, which can lead to different RBSK iterations.
  To obtain the optimal RBSK iteration, one needs to solve the following optimization
  problem
  \begin{eqnarray}
   \nonumber 
    [\bbb{P}_{\txt{opt}},S_{\txt{opt}}] &=& \arg\min_{\bbb{P}}\left(\min_{S\in\mathbb{DR}^{m\times m}} U(\bbb{P},S)\right) \\
    \label{upperBound_opt2}
    &=& \arg\min_{\bbb{P},S} U(\bbb{P},S),
  \end{eqnarray}
  or equivalently,
  \begin{eqnarray*}
    U(\bbb{P}_{\txt{opt}},S_{\txt{opt}}) &=& \min_{\bbb{P},S} U(\bbb{P},S).
  \end{eqnarray*}
  $U(\bbb{P}_{\txt{opt}},S_{\txt{opt}})$ is the sharpest upper bound of the convergence rate of the RBSK method that can be achieved for all possible joint distributions $\bbb{P}$. The matrix $S_{\txt{opt}}$ may be different
  from the matrix $\widehat{S}_{\txt{opt}}$. Obviously, the optimization problem {\rm(\ref{upperBound_opt2})} is
  very difficult. In order to simplify the resolution of {\rm(\ref{upperBound_opt2})}, an alternative plan is to
  fix the value of $S$ (e.g., let $\widehat{S}=\diag\{\frac{1}{\|A_1\|_2},\ldots,\frac{1}{\|A_m\|_2}\}$, then
  $\widehat{S}A$ has normalized rows; see such settings in {\rm\cite{NeedellTropp2014LAA})}, and solve the following
  optimization problem
  \begin{eqnarray}
  \label{upperBound_opt3}
    \widehat{\bbb{P}}_{\txt{opt}} &=& \arg\min_{\bbb{P}} U(\bbb{P},\widehat{S}).
  \end{eqnarray}
  Although {\rm(\ref{upperBound_opt3})} is a simplified version of {\rm(\ref{upperBound_opt2})}, it is still difficult to solve when the system size is large. In practice, various greedy and adaptive sampling rules have been designed to construct probability distributions that lead to improved convergence behavior by adjusting the sampling distribution dynamically according to current residual information. However, these strategies typically require substantial additional computation and data movement at each iteration in order to evaluate residual-related quantities or other metrics. An alternative idea, which we only briefly outline here, is to employ a learning-based approach: instead of recomputing probabilities adaptively during each iteration, one may train a model (e.g., neural networks) on a family of problem instances (for example, computed tomography problems) to produce a batch-sampling distribution $\bbb{P}$. Once such a distribution is obtained, it can be directly utilized in the RBSK iterations for solving the aforementioned problem instances. This strategy preserves low per-iteration costs of classical randomized methods in each RBSK iteration, while implicitly capturing problem-specific structures that can accelerate convergence. The effectiveness of this strategy has been preliminarily validated through experiments, although these results are not presented here. Developing and analyzing such learning-guided distributions lies beyond the scope of this paper, but it suggests an interesting direction for future research.

\end{remark}

\begin{remark}
\label{remark5.3}
  Due to the appearance of the scaling matrix $S$, the bounds in Theorem \ref{convRBSK} are independent on the magnitudes of the rows of the data matrix $A$ in {\rm(\ref{linear-system})}.
\end{remark}

The following corollary provides a relaxed convergence rate bound of the RBSK method.

\begin{corollary}
  \label{convRBSK_relaxed}
  Under the same conditions as Theorem \ref{convRBSK}, the mean squared error satisfies
  \begin{eqnarray}
    \label{thm2Result}
    \E\left(\|x^{(k+1)}-x_{\star}\|_2^2\right)
    &\leq& \left[\min_{S\in\mathbb{DR}^{m\times m}}\left(1-\widehat{\xi}\right)\right]^k
    \left[\min_{S\in\mathbb{DR}^{m\times m}}\left(1-\widehat{\eta}\right)\right]\|x^{(0)}-x_{\star}\|_2^2,
  \end{eqnarray}
  in particular, for $k=0$, it satisfies
  \begin{eqnarray}
     \label{thm2Eq5}
    \E\left(\|x^{(1)}-x_{\star}\|_2^2\right) 
    &\leq& \left[\min_{S\in\mathbb{DR}^{m\times m}}\left(1-\widehat{\eta}\right)\right]\|x^{(0)}-x_{\star}\|_2^2,
  \end{eqnarray}
  and, for $k>0$, it satisfies
  \begin{eqnarray}
    \label{thm2Eq6}
    \E\left(\|x^{(k+1)}-x_{\star}\|_2^2\right)
    &\leq& \left[\min_{S\in\mathbb{DR}^{m\times m}}\left(1-\widehat{\xi}\right)\right]\E\left(\|x^{(k)}-x_{\star}\|_2^2\right),
  \end{eqnarray}
  where
  \begin{eqnarray*}
    \widehat{\eta} &=&\frac{1}{\beta^S}\lambda_{\min}\left(A^{\T}S^{\T}\widehat{P}SA\right),
  \end{eqnarray*}
  and
  \begin{eqnarray*}
    \widehat{\xi} &=& \frac{1}{\beta^S}\min_{\tau\sim\bbb{P}} \lambda_{\min}\left(W_{\tau}^{\T}A^{\T} S^{\T}\widehat{P}SA W_{\tau}\right)
  \end{eqnarray*}
  with $\beta^S=\max_{\tau\sim\bbb{P}} \|S_{\tau}A_{\tau}\|_2^2$ and
  \begin{eqnarray}
  \label{thm2Eq0}
    \widehat{P} &=& \begin{bmatrix}
                   \sum_{\tau_{\txttiny{E$(1)$}}\in\range(\tau)} \pr\left(\tau_{\txttiny{E$(1)$}}\right)   &   &   \\
                      & \ddots &   \\
                      &   & \sum_{\tau_{\txttiny{E$(m)$}}\in\range(\tau)} \pr\left(\tau_{\txttiny{E$(m)$}}\right)
                  \end{bmatrix},
  \end{eqnarray}
  here, $S\in\mathbb{DR}^{m\times m}$ is any prescribed nonsingular diagonal matrix,
  and $\tau_{\txttiny{E$(j)$}}\in\range(\tau)$ are effective sets
  that include the row index $j$ (i.e., $j\in\tau_{\txttiny{E$(j)$}}$), for $1\leq j\leq m$.
\end{corollary}
{\em Proof.}
According to the definition of $B_{S;i}$ in (\ref{lem5Eq3}), it reads that
\begin{eqnarray*}
  B_{S;i}^{-1} &\succeq& \frac{1}{\beta^S} I,
\end{eqnarray*}
then, it follows that
\begin{eqnarray}
 \label{thm2Eq1a}
  D^2 &=& S^{\T}\left(\sum_{i=1}^{q} B_{S;i}^{-1} P_i\right)S \\ \nonumber
  &\succeq& \frac{1}{\beta^S} S^{\T} \left(\sum_{i=1}^{q} P_i\right)S \\ \label{thm2Eq1}
  &=& \frac{1}{\beta^S} S^{\T}\widetilde{P} S,
\end{eqnarray}
where
\begin{eqnarray*}
  \widetilde{P} &=& \begin{bmatrix}
                      \sum_{i=1}^{q} p_{i1} &   &   \\
                        & \ddots &   \\
                        &   & \sum_{i=1}^{q} p_{im}
                    \end{bmatrix}.
\end{eqnarray*}
According to Lemma \ref{lem8}, it holds that
\begin{eqnarray*}
  \sum_{i=1}^{q} p_{ij} &\geq& \sum_{\tau_{\txttiny{E$(j)$}}\in\range(\tau)} \pr\left(\tau_{\txttiny{E$(j)$}}\right),\
    \mbox{for }\ j=1,\ldots,m,
\end{eqnarray*}
which leads to a fact
\begin{eqnarray}
\label{thm2Eq2}
  \widetilde{P} &\succeq& \widehat{P}.
\end{eqnarray}
The relations (\ref{thm2Eq1}) and (\ref{thm2Eq2}) imply the following fact
\begin{eqnarray*}
  D^2 &\succeq& \frac{1}{\beta^S} S^{\T} \widehat{P}S.
\end{eqnarray*}
Due to the above relation, it follows that
\begin{eqnarray}
 \nonumber 
  1-\eta &=& 1-\lambda_{\min}\left(A^{\T}D^2A\right) \\ \label{thm2Eq3}
  &\leq& 1-\frac{1}{\beta^S}\lambda_{\min}\left(A^{\T}S^{\T}\widehat{P}SA\right)
\end{eqnarray}
and
\begin{eqnarray}
 \nonumber 
  1-\xi &=& 1-\min_{\tau\sim\bbb{P}} \lambda_{\min}\left(W_{\tau}^{\T}A^{\T} D^2AW_{\tau}\right) \\ \label{thm2Eq4}
  &\leq& 1- \frac{1}{\beta^S}\min_{\tau\sim\bbb{P}} \lambda_{\min}\left(W_{\tau}^{\T}A^{\T} S^{\T}\widehat{P}SAW_{\tau}\right),
\end{eqnarray}
where $\eta$ and $\xi$ are defined in Theorem \ref{convRBSK}.
For $k=0$, the inequality (\ref{thm2Eq3}) and the estimate (\ref{thm1Eq2}) lead to the relaxed estimate (\ref{thm2Eq5}).
For $k>0$, the inequality (\ref{thm2Eq4}) and the estimate (\ref{thm1Eq5}) lead to the relaxed estimate (\ref{thm2Eq6}).
Combining the estimates (\ref{thm2Eq5}) and (\ref{thm2Eq6}), and through recursion with respect to $k$, it can be concluded that
the mean squared error of the iterative sequence $\{x^{(k)}\}_{k=0}^{+\infty}$ obeys (\ref{thm2Result}).
$\hfill\square$

\begin{remark} \label{remark:relaxOne}
  In the cases of Example \ref{exampRowPaving} (\textsc{row paving batch-sampling}),
  Example \ref{exampUniform} (\textsc{uniform batch-sampling}), and Example \ref{exampNonUnique} (\textsc{non-unique batch-sampling}), the constant $\beta^S$ involved in the relaxed bound in Corollary \ref{convRBSK_relaxed} reads that
  \begin{eqnarray}
  \label{betaS_value}
    \beta^S &=& \left\{\begin{array}{rl}
                         \max_{\tau_{\txttiny{E}}\in\mathcal{T}}\|S_{\tau_\txttiny{E}}A_{\tau_\txttiny{E}}\|_2^2, & \txt{Example \ref{exampRowPaving}}, \\
                         \max_{\tau_{\txttiny{E}}\in\mathcal{F}}\|S_{\tau_\txttiny{E}}A_{\tau_\txttiny{E}}\|_2^2, & \txt{Example \ref{exampUniform}}, \\
                         \max_{\tau_{\txttiny{E}}\in\mathcal{G}}\|S_{\tau_\txttiny{E}}A_{\tau_\txttiny{E}}\|_2^2, & \txt{Example \ref{exampNonUnique}},
                       \end{array}
    \right.
  \end{eqnarray}
  where $\mathcal{T}$, $\mathcal{F}$, and $\mathcal{G}$ are the ranges of the batch-sampling $\tau$ defined in Example \ref{exampRowPaving}, Example \ref{exampUniform}, and Example \ref{exampNonUnique}. In addition, when the effective sets $\tau_{\txttiny{E}}$ are sampled uniformly at random in Examples
  \ref{exampRowPaving}, \ref{exampUniform}, and \ref{exampNonUnique}, the matrix $\widehat{P}$ involved in the relaxed bound in Corollary \ref{convRBSK_relaxed} reads that
  \begin{eqnarray}
  \label{Phat_value}
    \widehat{P} &=& \left\{\begin{array}{rl}
                             \frac{1}{L}\cdot I, & \txt{Example \ref{exampRowPaving}}, \\
                             \frac{\tbinom{m-1}{q-1}}{\tbinom{m}{q}}\cdot I, & \txt{Example \ref{exampUniform}}, \\
                             \frac{\sum_{t=1}^{q}\tbinom{m-1}{t-1}}{\sum_{t=1}^{q}\tbinom{m}{t}}\cdot I, & \txt{Example \ref{exampNonUnique}}.
                           \end{array}
    \right.
  \end{eqnarray}

  For the case of Example \ref{exampRowPaving} (\textsc{row paving batch-sampling}), let the matrix $S$ be fixed as $S=I$, together with the 1st lines of {\rm(\ref{betaS_value})} and {\rm(\ref{Phat_value})}, when $k=0$,
  Corollary \ref{convRBSK_relaxed} leads to a relaxed bound as
  \begin{eqnarray}
  \label{rlx_upperBound_1}
    1-\frac{1}{\beta^IL} \lambda_{\min}\left(A^{\T}A\right),
  \end{eqnarray}
  where $\beta^I=\max_{\tau_{\txttiny{E}}\in\mathcal{T}}\|A_{\tau_{\txttiny{E}}}\|_2^2$ can serve as a value of $\beta_{\txt{\tiny Up}}$ in {\rm(\ref{def:rowPaving})}.
  When $k>0$, Corollary \ref{convRBSK_relaxed} leads to a relaxed bound as
  \begin{eqnarray}
  \label{rlx_upperBound_2}
    1-\frac{1}{\beta^IL} \min_{\tau_{\txttiny{E}}\in\mathcal{T}}\lambda_{\min}\left(W_{\tau_{\txttiny{E}}}^{\T}A^{\T}AW_{\tau_{\txttiny{E}}}\right).
  \end{eqnarray}

  The relaxed bound {\rm(\ref{rlx_upperBound_1})} for $k=0$   recovers the convergence rate bound of the RBK method provided by Needell in {\rm\cite{NeedellTropp2014LAA}}. However, the argument suggested by
  Needell in {\rm\cite{NeedellTropp2014LAA}} says that the upper bound {\rm(\ref{rlx_upperBound_1})} is valid for the RBK method only when the data matrix $A$ has normalized rows. Our argument shows that the bound {\rm(\ref{rlx_upperBound_1})} works for all cases of
  the data matrix.

  The relaxed bound {\rm(\ref{rlx_upperBound_2})} for $k>0$ is sharper than {\rm(\ref{rlx_upperBound_1})} due to the fact that
  \begin{eqnarray*}
    \lambda_{\min}\left(W_{\tau_{\txttiny{E}}}^{\T}A^{\T}AW_{\tau_{\txttiny{E}}}\right) &\geq&
    \lambda_{\min}\left(A^{\T}A\right).
  \end{eqnarray*}
  In addition, thanks to the arbitrariness of the scaling matrix $S$, the bounds {\rm(\ref{thm2Eq5})} and {\rm(\ref{thm2Eq6})} provided by Corollary \ref{convRBSK_relaxed} may be even sharper than {\rm(\ref{rlx_upperBound_1})} and {\rm(\ref{rlx_upperBound_2})}.
  Thus, the bounds provided by Theorem \ref{convRBSK} and Corollary \ref{convRBSK_relaxed} are all sharper than the bound proved by Needell in {\rm\cite{NeedellTropp2014LAA}}.
\end{remark}

According to the proof of Corollary \ref{convRBSK_relaxed}, the upper bounds (\ref{thm2Eq5}) and (\ref{thm2Eq6}) are
obtained by relaxations of the matrix $B_{S;i}^{-1}$ in
(\ref{thm2Eq1a}) and the matrix $\widetilde{P}$ in (\ref{thm2Eq1}). As a matter of fact, when the case of \textsc{row paving batch-sampling} is considered
in the RBSK method (equivalent to the RBK method), the matrix $B_{S;i}^{-1}$ can be exactly described.
Therefore, a new bound between the bounds of Theorems \ref{convRBSK} and Corollary \ref{convRBSK_relaxed} can be derived
by only making the relaxation of $\widetilde{P}$, which leads to the following Theorem \ref{convRBSK_relaxed_pave}.

\begin{theorem}
  \label{convRBSK_relaxed_pave}
  Under the same conditions as Theorem \ref{convRBSK}, let the marginal distributions of $\tau_i$ satisfy $p_{ij}>0$. Let $\tau\sim\bbb{P}$ refers to a
  \textsc{row paving batch-sampling}, and the range of $\tau$ (a partition $\mathcal{T}=\{T_1,\ldots,T_{L}\}$) defines a row paving ($L,\beta_{\txt{\tiny Low}},\beta_{\txt{\tiny Up}}$) of $A$. Then the mean squared error satisfies
  \begin{eqnarray}
    \label{thm3Result}
    \E\left(\|x^{(k+1)}-x_{\star}\|_2^2\right)
    &\leq& \left[\min_{S\in\mathbb{DR}^{m\times m}}\left(1-\widetilde{\xi}\right)\right]^k
    \left[\min_{S\in\mathbb{DR}^{m\times m}}\left(1-\widetilde{\eta}\right)\right]\|x^{(0)}-x_{\star}\|_2^2,
  \end{eqnarray}
  in particular, for $k=0$, it satisfies
  \begin{eqnarray}
     \label{thm3Eq3}
    \E\left(\|x^{(1)}-x_{\star}\|_2^2\right) 
    &\leq& \left[\min_{S\in\mathbb{DR}^{m\times m}}\left(1-\widetilde{\eta}\right)\right]\|x^{(0)}-x_{\star}\|_2^2,
  \end{eqnarray}
  and, for $k>0$, it satisfies
  \begin{eqnarray}
    \label{thm3Eq4}
    \E\left(\|x^{(k+1)}-x_{\star}\|_2^2\right)
    &\leq& \left[\min_{S\in\mathbb{DR}^{m\times m}}\left(1-\widetilde{\xi}\right)\right]\E\left(\|x^{(k)}-x_{\star}\|_2^2\right),
  \end{eqnarray}
  where
  \begin{eqnarray*}
    \widetilde{\eta} &=&\lambda_{\min}\left(A^{\T}S^{\T}B_S^{-1} \widehat{P}SA\right),
  \end{eqnarray*}
  and
  \begin{eqnarray*}
    \widetilde{\xi} &=& \min_{\tau\sim\bbb{P}} \lambda_{\min}\left(W_{\tau}^{\T}A^{\T} S^{\T} B_S^{-1} \widehat{P}SA W_{\tau}\right)
  \end{eqnarray*}
  with $\widehat{P}$ defined in {\rm(\ref{thm2Eq0})}, in addition, when an indicator $\mathcal{I}(\cdot)$
  is given by $\mathcal{I}(j)=t$, if $j\in T_t\subset\mathcal{T}$, for $j=1,\ldots,m$,
  the matrix $B_S$ is of the form
  \begin{eqnarray*}
    B_S &=& \begin{bmatrix}
              \|S_{T_{\mathcal{I}(1)}}A_{T_{\mathcal{I}(1)}}\|_2^2 &   &   \\
                & \ddots &   \\
                &   & \|S_{T_{\mathcal{I}(m)}}A_{T_{\mathcal{I}(m)}}\|_2^2
            \end{bmatrix},
  \end{eqnarray*}
  here, $S\in\mathbb{DR}^{m\times m}$ is any prescribed nonsingular diagonal matrix.
\end{theorem}
{\em Proof.}
Since the marginal distributions $p_{ij}$ are positive,
the quantities $\beta_{i,j}^S$ introduced in Lemma \ref{lem6} are of the form
\begin{eqnarray*}
  \beta_{i,j}^S &=& \|S_{T_{\mathcal{I}(j)}}A_{T_{\mathcal{I}(j)}}\|_2^2,\ \forall j=1,\ldots,m,
\end{eqnarray*}
for $1\leq i\leq q$. Obviously, $\beta_{i,j}^S$ is independent of $i$, together with
the definition of $B_{S;i}$ in (\ref{lem5Eq3}), it leads to a fact that
\begin{eqnarray*}
  B_{S;i} &=& B_S,\ \forall 1\leq i\leq q.
\end{eqnarray*}
According to the proofs of Theorems \ref{convRBSK} and Corollary \ref{convRBSK_relaxed}, it follows that
\begin{eqnarray*}
  D^2 &=& S^{\T}\left(\sum_{i=1}^{q} B_{S;i}^{-1} P_i\right)S \\
  &=& S^{\T}\left(B_S^{-1} \sum_{i=1}^{q} P_i\right)S \\
  &\succeq& S^{\T}B_S^{-1} \widehat{P}S.
\end{eqnarray*}
Due to the above relation, it holds that
\begin{eqnarray}
 \nonumber 
  1-\eta &=& 1-\lambda_{\min}\left(A^{\T}D^2A\right) \\ \label{thm3Eq1}
  &\leq& 1-\lambda_{\min}\left(A^{\T}S^{\T}B_S^{-1} \widehat{P}SA\right)
\end{eqnarray}
and
\begin{eqnarray}
 \nonumber 
  1-\xi &=& 1-\min_{\tau\sim\bbb{P}} \lambda_{\min}\left(W_{\tau}^{\T}A^{\T} D^2AW_{\tau}\right) \\ \label{thm3Eq2}
  &\leq& 1- \min_{\tau\sim\bbb{P}} \lambda_{\min}\left(W_{\tau}^{\T}A^{\T} S^{\T} B_S^{-1} \widehat{P}SAW_{\tau}\right),
\end{eqnarray}
where $\eta$ and $\xi$ are defined in Theorem \ref{convRBSK}.
For $k=0$, the inequality (\ref{thm3Eq1}) and the estimate (\ref{thm1Eq2}) lead to the new estimate (\ref{thm3Eq3}).
For $k>0$, the inequality (\ref{thm3Eq2}) and the estimate (\ref{thm1Eq5}) lead to the new estimate (\ref{thm3Eq4}).
Combining the new estimates (\ref{thm3Eq3}) and (\ref{thm3Eq4}), and through recursion with respect to $k$, it can be concluded that
the mean squared error of the iterative sequence $\{x^{(k)}\}_{k=0}^{+\infty}$ obeys (\ref{thm3Result}).
$\hfill\square$

\begin{remark} \label{remark:relaxTwo}
  Obviously, in Theorem \ref{convRBSK_relaxed_pave}, due to the fact that the batch-sampling $\tau$ refer to the case of \textsc{row paving batch-sampling}, the matrix $B_S$ satisfies
  \begin{eqnarray*}
    B_S^{-1} &\succeq& \frac{1}{\beta^S}I,
  \end{eqnarray*}
  where
  \begin{eqnarray*}
    \beta^S &=& \max_{\tau\sim\bbb{P}} \|S_{\tau} A_{\tau}\|_2^2 \\
    &=& \max_{1\le j\le m} \|S_{T_{\mathcal{I}(j)}}A_{T_{\mathcal{I}(j)}}\|_2^2 \\
    &=& \max_{1\le t\le L} \|S_{T_t}A_{T_t}\|_2^2. \\
  \end{eqnarray*}
  Thus, the bounds in Theorem \ref{convRBSK_relaxed_pave} are sharper than those in Corollary \ref{convRBSK_relaxed} and Needell's bounds in {\rm\cite{NeedellTropp2014LAA}}.
\end{remark}

Based on Theorem \ref{convRBSK}, we can derive an even sharper convergence rate bound of the RBSK method via concentration inequalities, under slightly stronger assumptions. The derivation here is motivated by the idea behind the law of large numbers. Specifically, for a sufficiently large sample size, the sample mean can closely approximate the expectation of the random variable. As a result, the refined theoretical rate bound is expressed relying on the expectation rather than the sample mean.
\begin{theorem}
	\label{convRBSK_sharper}
	Under the same conditions as Theorem \ref{convRBSK},
	for $k>0$, consider the random variable
	$\xi_{\tau^{(k-1)}}=\lambda_{\min}\left(W_{\tau^{(k-1)}}^{\T}A^{\T} D^2AW_{\tau^{(k-1)}}\right).$ Let $\left\{\xi_{\tau^{(k-1,j)}}\right\}_{j=1}^{\ell}$ be a sample sequence of $\xi_{\tau^{(k-1)}}$, obtained from independent realizations $\left\{\tau^{(k-1, j)}\right\}_{j=1}^{\ell}$ of $\tau^{(k-1)}$.
	Similarly, define $\left\{\|x^{(k, j)}-x_{\star}\|_2\right\}_{j=1}^{\ell}$ as the corresponding sample sequence of the error norms $\|x^{(k)}-x_{\star}\|_2$, where $x^{(k, j)}$ denotes the $k$th iterate generated from the $j$th realization. Let $a=\lambda_{\min}(A^{\T}D^2A)$, $b=\lambda_{\max}(A^{\T}D^2A)$. For any $\epsilon>0, \delta \in (0,1)$, if the sample size $\ell$ reads that
	\begin{eqnarray*}
		\ell &\geq& \frac{(b - a)^2}{2 \epsilon^2} \log\left(\frac{2}{\delta}\right),
	\end{eqnarray*}
	and the sample covariance satisfies $\cov\left(\xi_{\tau^{(k-1)}},\|x^{(k)}-x_{\star}\|^2_2\right)\geq 0$, then, with probability at least $1-\delta$, the mean squared error satisfies
	\begin{eqnarray}
		\label{thm4Eq2}
		\E\left(\|x^{(k+1)}-x_{\star}\|_2^2\right)
		&\leq& \left[\min_{S\in\mathbb{DR}^{m\times m}}\left(1+\epsilon-\E\left(\xi_{\tau} \right)\right)   \right]\E\left(\|x^{(k)}-x_{\star}\|_2^2\right),
	\end{eqnarray}	
	moreover, when $\epsilon,\delta\to 0$, it follows with probability one that
	\begin{eqnarray}
		\label{thm4Eqfinal}
		\E\left(\|x^{(k+1)}-x_{\star}\|_2^2\right)
		&\leq&
		\left[\min_{S\in\mathbb{DR}^{m\times m}}\left(1-\E\left(\xi_{\tau}\right)\right)\right]
		\E\left(\|x^{(k)}-x_{\star}\|_2^2\right),
	\end{eqnarray}
	where $\xi_{\tau}=\lambda_{\min}\left(W_{\tau}^{\T}A^{\T} D^2AW_{\tau}\right)$, and $D$ is the diagonal matrix defined in \eqref{D_formula}.	
\end{theorem}
{\em Proof.}
According to the inequality (\ref{thm1Eq3}) in the proof of Theorem \ref{convRBSK}, it follows that
\begin{eqnarray}
	\label{thm4Eq3}
	\E_k\left(\|x^{(k+1,j)}-x_{\star}\|_2^2\right)
	&\leq& \left(1-\xi_{\tau^{(k-1,j)}}\right)\|x^{(k,j)}-x_{\star}\|_2^2,\ \forall 1\leq j\leq \ell.
\end{eqnarray}
with
\begin{eqnarray*}
	\xi_{\tau^{(k-1,j)}} &=& \lambda_{\min}\left(W_{\tau^{(k-1,j)}}^{\T}A^{\T} D^2AW_{\tau^{(k-1,j)}}\right).
\end{eqnarray*}
Taking the sample mean on both sides of inequality (\ref{thm4Eq3}), we obtain that
\begin{eqnarray}
	\label{thm4Eq4}
	\frac{1}{\ell}\sum_{j=1}^{\ell}\E_k\left(\|x^{(k+1,j)}-x_{\star}\|_2^2\right)
	&\leq& \frac{1}{\ell}\sum_{j=1}^{\ell}\left[\left(1-\xi_{\tau^{(k-1,j)}}\right)\|x^{(k,j)}-x_{\star}\|_2^2\right].
\end{eqnarray}
According to Lemma \ref{lem5}, since the sample covariance satisfies $\cov\left(\xi_{\tau^{(k-1)}},\|x^{(k)}-x_{\star}\|^2_2\right)\geq0$, then
\begin{eqnarray}
	\label{thm4Eq5}
 	\frac{1}{\ell}\sum_{j=1}^{\ell}\xi_{\tau^{(k-1,j)}}\|x^{(k,j)}-x_{\star}\|_2^2 \geq \left(\frac{1}{\ell}\sum_{j=1}^{\ell}\xi_{\tau^{(k-1,j)}}\right)\left(\frac{1}{\ell}\sum_{j=1}^{\ell}\|x^{(k,j)}-x_{\star}\|_2^2\right).
\end{eqnarray}
Combining inequalities (\ref{thm4Eq4}) and (\ref{thm4Eq5}), it follows that
\begin{eqnarray}
	\label{thm4Eq6}
	\frac{1}{\ell}\sum_{j=1}^{\ell}\E_k\left(\|x^{(k+1,j)}-x_{\star}\|_2^2\right)
	&\leq& \left(1-\frac{1}{\ell}\sum_{j=1}^{\ell}\xi_{\tau^{(k-1,j)}}\right) \left(\frac{1}{\ell}\sum_{j=1}^{\ell}\|x^{(k,j)}-x_{\star}\|_2^2\right).
\end{eqnarray}
By the definition of $\xi_{\tau^{(k-1,j)}}$, it can be readily verified that
\begin{eqnarray*}
 	a \leq \xi_{\tau^{(k-1,j)}} \leq b.
\end{eqnarray*}
According to Corollary \ref{Hoeffding-co}, for any $\epsilon>0, \delta \in (0,1)$, since the sample size $\ell$ satisfies
\begin{eqnarray*}
	\ell &\geq& \frac{(b - a)^2}{2 \epsilon^2} \log\left(\frac{2}{\delta}\right),
\end{eqnarray*}
then it follows that
\begin{eqnarray*}
	\pr\left(\left|\frac{1}{\ell}\sum_{j=1}^{\ell}\xi_{\tau^{(k-1,j)}}-\E(\xi_{\tau^{(k-1)}})\right|\geq\epsilon\right)  &\leq& \delta. \\
\end{eqnarray*}
Consequently, this implies that,
\begin{eqnarray}
	\label{thm4Eq7}
	\pr\left(\frac{1}{\ell}\sum_{j=1}^{\ell}\xi_{\tau^{(k-1,j)}}-\E(\xi_{\tau^{(k-1)}})\geq -\epsilon \right)  &\geq& 1-\delta
\end{eqnarray}
Due to the relation (\ref{thm4Eq6}) and (\ref{thm4Eq7}), we obtain the inequality
\begin{eqnarray}
	\label{thm4Eq8}
	\frac{1}{\ell}\sum_{j=1}^{\ell}\E_k\left(\|x^{(k+1,j)}-x_{\star}\|_2^2\right)
	&\leq& \left(1+\epsilon-\E\left(\xi_{\tau^{(k-1)}} \right)\right) \left(\frac{1}{\ell}\sum_{j=1}^{\ell}\|x^{(k,j)}-x_{\star}\|_2^2\right).
\end{eqnarray}
which holds with probability at least $1-\delta$.
According to the notation in Section \ref{sec:prelimNota}, it reads that
\begin{align*}
	\E_k\left(\|x^{(k+1,j)}-x_{\star}\|_2^2\right)
	&=\E\left(\|x^{(k+1,j)}-x_{\star}\|_2^2\middle|\tau^{(0)},\ldots,\tau^{(k-2)},\tau^{(k-1,j)}\right),\\
	\|x^{(k,j)}-x_{\star}\|_2^2
	&=\E\left(\|x^{(k,j)}-x_{\star}\|_2^2\middle|\tau^{(0)},\ldots,\tau^{(k-2)},\tau^{(k-1,j)}\right),
\end{align*}
with $\tau^{(0)},\ldots,\tau^{(k-2)},\tau^{(k-1)},\tau^{(k-1,j)} \sim \bbb{P}$. Thus, the full expectations of the above errors hold that
\begin{align*}
	\E\left(\E_k\left(\|x^{(k+1,j)}-x_{\star}\|_2^2\right)\right)
	&=\E\left(\E\left(\|x^{(k+1,j)}-x_{\star}\|_2^2\middle|\tau^{(0)},\ldots,\tau^{(k-2)},\tau^{(k-1,j)}\right)\right)\\
	&=\E\left(\E\left(\|x^{(k+1)}-x_{\star}\|_2^2\middle|\tau^{(0)},\ldots,\tau^{(k-2)},\tau^{(k-1)}\right)\right)\\
	&=\E\left(\|x^{(k+1)}-x_{\star}\|_2^2\right),
\end{align*}
and
\begin{align*}
	\E\left(\|x^{(k,j)}-x_{\star}\|_2^2\right)
	&=\E\left(\E\left(\|x^{(k,j)}-x_{\star}\|_2^2\middle|\tau^{(0)},\ldots,\tau^{(k-2)},\tau^{(k-1,j)}\right)\right)\\
	&=\E\left(\E\left(\|x^{(k)}-x_{\star}\|_2^2\middle|\tau^{(0)},\ldots,\tau^{(k-2)},\tau^{(k-1)}\right)\right)\\
	&=\E\left(\|x^{(k)}-x_{\star}\|_2^2\right).
\end{align*}
Since $\tau$ and $\tau^{(k-1)}$ follow the same distribution, we can replace $\E(\xi_{\tau^{(k-1)}})$ by $\E(\xi_{\tau})$. Therefore, taking full expectation on both sides of \eqref{thm4Eq8} and minimizing with respect to $S$ yields the estimate \eqref{thm4Eq2}.

Moreover, let $\epsilon,\delta \to 0$, then $\ell \to +\infty$. By the strong law of large numbers, it holds almost surely that
\begin{align*}
	\lim_{\ell\to+\infty}\frac{1}{\ell}\sum_{j=1}^{\ell}\E\left(\|x^{(k+1,j)}-x_{\star}\|_2^2
	\middle|\tau^{(0)},\ldots,\tau^{(k-2)},\tau^{(k-1,j)}\right)
	&\overset{\mathrm{a.s.}}{=}\E\left(\|x^{(k+1)}-x_{\star}\|_2^2\middle|\tau^{(0)},\ldots,\tau^{(k-2)}\right),
\end{align*}
and
\begin{align*}
	\lim_{\ell\to+\infty}\frac{1}{\ell}\sum_{j=1}^{\ell}\E\left(\|x^{(k,j)}-x_{\star}\|_2^2\middle|\tau^{(0)},\ldots,\tau^{(k-2)},\tau^{(k-1,j)}\right)
	&\overset{\mathrm{a.s.}}{=}\E\left(\|x^{(k)}-x_{\star}\|_2^2\middle|\tau^{(0)},\ldots,\tau^{(k-2)}\right).
\end{align*}
Combining these relations with \eqref{thm4Eq8} and $\epsilon,\delta\to 0$, we obtain with probability one that
\begin{align*}
	\E\left(\|x^{(k+1)}-x_{\star}\|_2^2\middle|\tau^{(0)},\ldots,\tau^{(k-2)}\right)
	&\leq\left(1-\E(\xi_{\tau})\right)\E\left(\|x^{(k)}-x_{\star}\|_2^2\middle|\tau^{(0)},\ldots,\tau^{(k-2)}\right).
\end{align*}
Finally, by taking full expectation on both sides of the above estimate and minimizing with respect to $S$, we have \eqref{thm4Eqfinal} holds with probability one.
 $\hfill\square$

\begin{remark} \label{covariance:explain}
	We validate by numerical experiments that the covariance condition in Theorem~\ref{convRBSK_sharper}, that is,  $\cov\left(\xi_{\tau^{(k-1)}},\|x^{(k)}-x_{\star}\|_2^2\right)\geq 0,$
	occurs with high empirical frequency in RBSK iterations.
	To verify this, we design the following experiment. Since the initial guess $x^{(0)}\in\range(A^{\T})$, and each RBSK update is obtained by adding a vector in $\range(A^{\T})$, all iterates generated by RBSK remain in $\range(A^{\T})$. Moreover, the least-norm solution $x_{\star}$ also belongs to $\range(A^{\T})$. Therefore, to simulate arbitrary iterates arising in the RBSK update, we sample test points from a neighborhood of $x_{\star}$ within $\range(A^{\T})$. Specifically, for each test matrix, we randomly generate $100$ test points of the form $x_{\star}+\alpha \nu$, where $\nu\in \range(A^{\T})$ is a random unit vector and $\alpha \in [0, R]$ is a random scalar with $R = 10$. Thus, all test points lie in the ball centered at $x_{\star}$ with radius $R$. For each sampled test point $x^{(k-1)}$, we then draw $100$ independent batch samples $\left\{\tau^{(k-1, j)}\right\}_{j=1}^{100}$, compute the corresponding values of $\xi_{\tau^{(k-1,j)}}$, perform one RBSK update for each $j$, and compute the values of $\|x^{(k,j)}-x_{\star}\|_2^2$. The sample covariance between the two sample sequences is then recorded for each sampled test point. The experiment is carried out for two Gaussian random matrices of sizes $1000\times 100$ and $100\times 1000$, respectively, as shown in FIG.~\ref{fig:covariance}, with batch size $q = 10$. In both cases, the sample covariance is positive for the majority of the sampled test points, which indicates that the condition $\cov\left(\xi_{\tau^{(k-1)}},\|x^{(k)}-x_{\star}\|^2_2\right)\geq 0$ is satisfied with high empirical frequency. The experiment provides empirical support for the reasonableness of the covariance condition and hence for the applicability of Theorem~\ref{convRBSK_sharper}.
	\begin{figure}[H]
		\centering
		\subfigure[]{\includegraphics[width=0.45\linewidth]{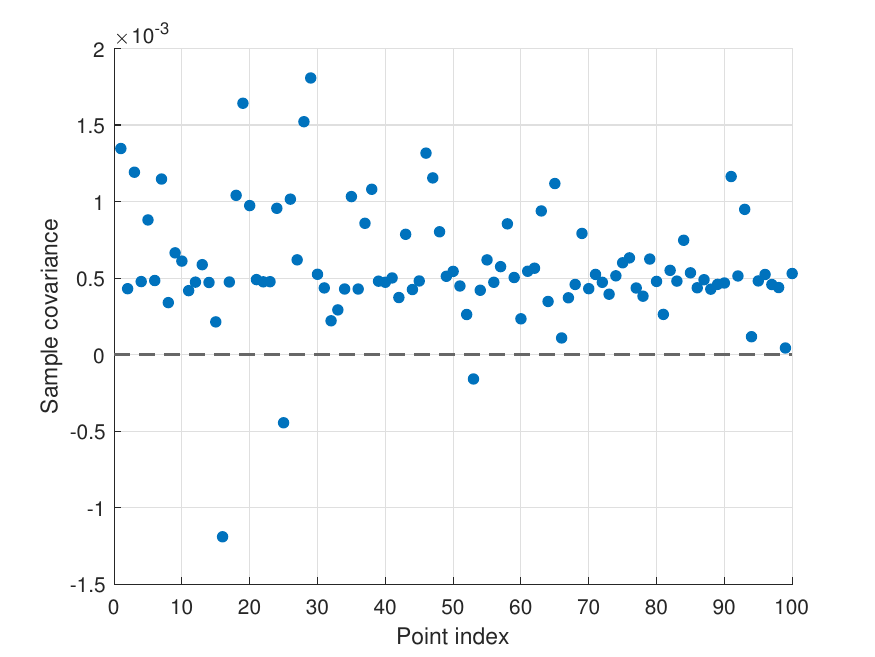}}
		\subfigure[]{\includegraphics[width=0.45\linewidth]{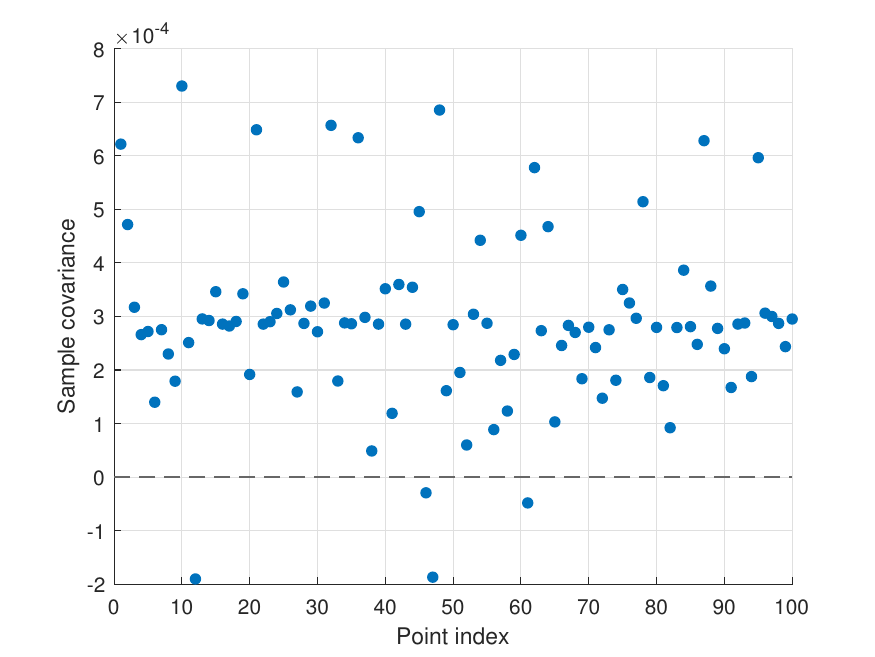}}
		\caption{Sample covariance $\cov\left(\xi_{\tau^{(k-1)}},\|x^{(k)}-x_{\star}\|_2^2\right)$ at sampled test points in the neighborhood of $x_{\star}$ with radius $R=10$. Subfigures (a) and (b) correspond to matrix sizes $100\times1000$ and $1000\times100$, with batch-size $q=10$.}
		\label{fig:covariance}
	\end{figure}	
\end{remark}

\section{Experimental Results}
\label{sec:expResult}
In this section, we report the numerical results that compare the theoretical convergence rate bound derived in the new batch-sampling framework of RBSK with two existing results: the classical bound established by Needell and Tropp in \cite{NeedellTropp2014LAA} (row-paving framework), hereafter referred to as ND14, and the bound proposed by Gower et al. in \cite{gower2021adaptive} (sketch-and-project framework), hereafter referred to as GM21. These comparisons provide numerical evidence supporting the effectiveness of our theoretical results.

We will perform batch-sampling on the data matrix in the experiments in a manner compatible with ND14 and GM21. The reasons are as follows. The RBSK framework only requires knowledge of the batch-sampling distribution, whereas ND14 and GM21 require the explicit pre-construction of all possible sub-blocks of the data matrix.
In particular, ND14 further requires the rows to be partitioned into mutually disjoint blocks, so that each block has distinct rows.
Due to the generality of the RBSK framework, any sampling rule defined in ND14 or GM21 can be directly mapped to a sampling distribution in RBSK, though the converse is not generally true. Additionally, since adaptive batch-sampling is not within the scope of this paper, we only compare the bounds derived from the RBSK framework with ND14 and non-adaptive GM21.

The experiments focus on two types of data matrices: those with multi-scale structure and those that are ill-conditioned. For each type, we construct test instances using both synthetic and real-world data matrices. The detailed experimental setup and the exact formulas used for the test bounds are given in the following subsection.

\subsection{Experimental setup}
For synthetic data matrices, we use the \texttt{randn} function in MATLAB to create matrices whose entries follow the standard Gaussian distribution. Real-world data matrices are primarily obtained from the SuiteSparse Matrix Collection \cite{davis2011university}, which includes sparse matrices arising from various practical applications. To ensure the consistency of the test linear system, we generate a true solution by setting $x_{\star}=\texttt{randn(n,1)}$ and define the right-hand side as $b=Ax_{\star}$. To ensure compatibility with the results of ND14 and GM21, we employ a row-paving setting for the sub-data blocks instead of the generalized setting under the RBSK framework. Specifically, the data matrix $A \in \mathbb{R}^{m \times n}$ is partitioned naturally by rows into $L=m/q$ blocks of size $q \times n$, denoted as $A=[A_{\tau(1)},A_{\tau(2)},\ldots,A_{\tau(L)}]^{\T}$. For the batch-sampling strategy, we follow the recommendation in \cite{gower2015randomized,gower2021adaptive} and sample sub-data blocks according to a static probability distribution $p_i \sim\left\| A_{\tau^{(i)}} \right\|^2$. The initial guess of the iteration is $x^{(0)}=0$ in all tests. We define the relative solution error (RSE) at the $k$th iteration as
\begin{eqnarray*} 					  	     							\mathrm{RSE}=\frac{\left\|x^{(k)}-A^{\dagger}b\right\|_2^2}{\left\|x^{(0)}-A^{\dagger} b\right\|_2^2} .
\end{eqnarray*}
The algorithm terminates when the RSE drops below $10^{-8}$, or the number of iterations reaches $5000$.
Under the above settings, the algorithms in the framework of RBSK, ND14 and GM21 reduce to the same block Kaczmarz method. In addition, Theorem \ref{convRBSK} reduces to the form of Theorem \ref{convRBSK_relaxed_pave}. Therefore, we focus on testing the bounds of the convergence rates provided by Theorems \ref{convRBSK_relaxed_pave} and \ref{convRBSK_sharper}, which are denoted as `Thm \ref{convRBSK_relaxed_pave} bound' and `Thm \ref{convRBSK_sharper} bound', respectively.

We now present the explicit formulas for each test bound, and describe how the required quantities are computed in practice under the above row-paving setting.

\noindent\textbf{ND14.} ND14 has already been discussed in Section \ref{sec:convAnal}, and its original form is given in \eqref{rlx_upperBound_1}, which was established in \cite{NeedellTropp2014LAA} under the assumption that the working block is sampled uniformly at random. We have already demonstrated in Remarks \ref{remark:relaxOne} and \ref{remark:relaxTwo} that
all theoretical upper bounds derived within the RBSK framework are sharper than ND14. For the row-paving setting, where the block partition $\mathcal{T}=\{\tau(1),\tau(2),\ldots,\tau(L)\}$ is fixed and the block $\tau(i)$ is selected with probability $p_i$ at each iteration, following the derivation in Lemma 2.2 of \cite{NeedellTropp2014LAA}, it can be readily derived that the ND14 expression plotted in the experiments corresponds to
$$
1-\rho_{\text{\tiny ND14}}= 1-\lambda_{\min}\left(A^{\top} \frac{1}{\beta^I}\widehat{P} A\right),
$$
with $\beta^I=\max_i \|A_{\tau(i)}\|_2^2$, and $\widehat{P}$ will be given later. In particular, when $p_i=\frac{1}{L}$, $\forall i=1,\ldots,L$, this expression reduces to \eqref{rlx_upperBound_1}.

\noindent\textbf{GM21.} The expression for GM21 plotted in the figures is derived from the convergence rate bound of Gower et al. \cite[Section 7.2]{gower2021adaptive} specialized to the present settings. Specifically, GM21 can be written in the equivalent form
$$
1-\rho_{\text{\tiny GM21}}=1-\lambda_{\min}\left(A^{\top} \widehat{B}^{\dagger}\widehat{P} A\right)
$$
with $\widehat{B}=\text{blkdiag}(A_{\tau(1)}A_{\tau(1)}^{\top}, \ldots, A_{\tau(L)}A_{\tau(L)}^{\top})$.

\begin{table}[htbp]
	\renewcommand{\arraystretch}{1.4}
	\caption{Exact formulas of the test bounds in row-paving setting.}\label{tab:formulas}
	\centering
	\begin{tabular}{|c|c|c|c|}\hline
		Theoretical bound& Exact formula & Parameter specification &Reference\\\hline
		ND14&$1-\rho_{\text{\tiny ND14}}$&$\rho_{\text{\tiny ND14}}= \lambda_{\min}(A^{\top} \frac{1}{\beta^I}\widehat{P} A)$&\cite{NeedellTropp2014LAA}\\\hline
		
		GM21&$1-\rho_{\text{\tiny GM21}}$&$\rho_{\text{\tiny GM21}}=\lambda_{\min}(A^{\top} \widehat{B}^{\dagger}\widehat{P} A)$&\cite{gower2021adaptive} \\\hline
		
		Thm \ref{convRBSK_relaxed_pave} &$\displaystyle \min_{S\in\mathbb{DR}^{m\times m}}(1-\widetilde{\xi})$&$\displaystyle \widetilde{\xi}=\min_{\tau\sim\bbb{P}}\lambda_{\min}(W_{\tau}^{\T}A^{\T} S^2 B_S^{-1} \widehat{P}A W_{\tau})$&\eqref{thm3Eq4}\\\hline
		
		Thm \ref{convRBSK_sharper} &$\displaystyle \min_{S\in\mathbb{DR}^{m\times m}}(1-\E\left(\xi_{\tau} \right))$&$\xi_{\tau}=\lambda_{\min}(W_{\tau}^{\T}A^{\T} D^2AW_{\tau})$, $\tau\sim\bbb{P}$&\eqref{thm4Eq2}\\\hline
	\end{tabular}
\end{table}

In all figures, the exact formulas of all test bounds are summarized in Table \ref{tab:formulas}. The GM21 bound refines the ND14 bound by replacing the matrix $\frac{1}{\beta^I}\widehat{P}$ with the matrix $\widehat{B}^{\dagger}\widehat{P}$. In the Thm 5.2 and Thm 5.3 bounds, the matrices $S^2 B_S^{-1}$ and $D^2$ serve analogous roles to $\widehat{B}^{\dagger}\widehat{P}$, while $W_{\tau}$ constrains a matrix to a subspace of $\reals^n$. The Thm 5.3 bound further yields a tighter non-worst-case bound by taking the expectation over $\xi_{\tau}$ and minimizing over $S$. The matrices $S$, $W_{\tau}$, $\widehat{P}$ and $P_i$ are constructed as follows.

\textbf{Construction of $S$.}\quad In the experiments, we construct a fixed nonsingular diagonal matrix $S\in\mathbb{DR}^{m\times m}$ when evaluating the `Thm \ref{convRBSK_relaxed_pave} bound' and `Thm \ref{convRBSK_sharper} bound' for simplicity (see Remark \ref{remark5.2} for the justification). Specifically, we set $S=\diag\{\frac{1}{\|A^{(1)}\|_2},\ldots,\frac{1}{\|A^{(m)}\|_2}\}$, where each diagonal entry is the reciprocal of the $\ell_2$-norm of the corresponding row of $A$. We emphasize that this is only a naive choice of $S$, while a sharper bound could be obtained by optimizing over $S$.

\textbf{Construction of $W_{\tau}$.}\quad For each batch-sampling $\tau$, the matrix $W_{\tau}$ involved in the theoretical bounds can be computed as follows:
\begin{itemize}[itemsep=2pt, topsep=0pt, parsep=0pt, partopsep=0pt]
	\item compute an orthonormal basis matrix $Q\in\mathbb{R}^{n\times d}$ of $\range(A^{\top})$ via a thin QR factorization of $A^{\top}$;
	\item compute a matrix $Z_{\tau}$ whose columns form an orthonormal basis of $\zspace(A_{\tau}Q)$ making use of the MATLAB function \texttt{null}, where $\zspace(A_{\tau}Q)=\range\left((A_{\tau}Q)^{\top}\right)^{\perp}$;
	\item compute $W_{\tau}=QZ_{\tau}$, then the columns of $W_{\tau}$ form an orthonormal basis of $\range(A^{\top})\cap\zspace(A_{\tau})=\mathcal{W}_{\tau}$.
\end{itemize}

\textbf{Construction of $\widehat{P}$ and $P_i$.}\quad The remaining quantities appearing in Theorems~\ref{convRBSK_relaxed_pave} and \ref{convRBSK_sharper} can be computed explicitly as follows. The matrix $\widehat{P}$ appearing in Theorem~\ref{convRBSK_relaxed_pave} is given by
$$
\widehat{P}=\text{blkdiag}\{p_1I_{q\times q},\,p_2I_{q\times q},\,\ldots,\,p_LI_{q\times q}\},
$$
where $I_{q\times q}$ denotes the $q\times q$ identity matrix. Moreover, according to Definition \ref{def:batchSampling}, the matrices $P_1,\ldots,P_q$ in $D^2=S^{\top}\bigl(\sum_{i=1}^{q}B_{S;i}^{-1}P_i\bigr)S$ are determined by the marginal distribution of the batch-sampling vector. In the row-paving setting,
$P_i$ has the following block diagonal form
$$
P_i=\text{blkdiag}\{p_1E_{ii},\,p_2E_{ii},\,\ldots,\,p_LE_{ii}\},
$$
where $E_{ii}\in\mathbb{R}^{q\times q}$ denotes the matrix with 1 at the $(i,i)$ position and $0$ otherwise.

To visualize the empirical results, we adopt the black dashed line with square markers in the figures to represent the average convergence rate over 30 trials. The lightly shaded area signifies the range from the minimum to the maximum convergence rates, while the darker shaded region indicates the interquartile range (spanning from the 25th to the 75th percentile) of the empirical convergence rates. Moreover, the magenta solid line with circle markers represents the ND14 bound, the red solid line with square markers represents the GM21 bound, the cyan dashed line with cross markers represents the Thm \ref{convRBSK_relaxed_pave} bound, and the blue dashed line with asterisk markers represents the Thm \ref{convRBSK_sharper} bound.

\subsection{Multi-scale data matrices}
A data matrix $A$ is said to have a multi-scale structure if certain sub-data blocks differ from others in terms of magnitude. For simplicity, we construct matrices with a two-scale structure for testing purposes. Specifically, a sub-data block of $A$ is scaled by a factor $\alpha=0.2$, while leaving the other sub-data blocks unchanged. This operation introduces a clear magnitude disparity from one block to the others.

\begin{figure}[htbp]
	\centering
	\subfigure[]{\includegraphics[width=0.32\linewidth]{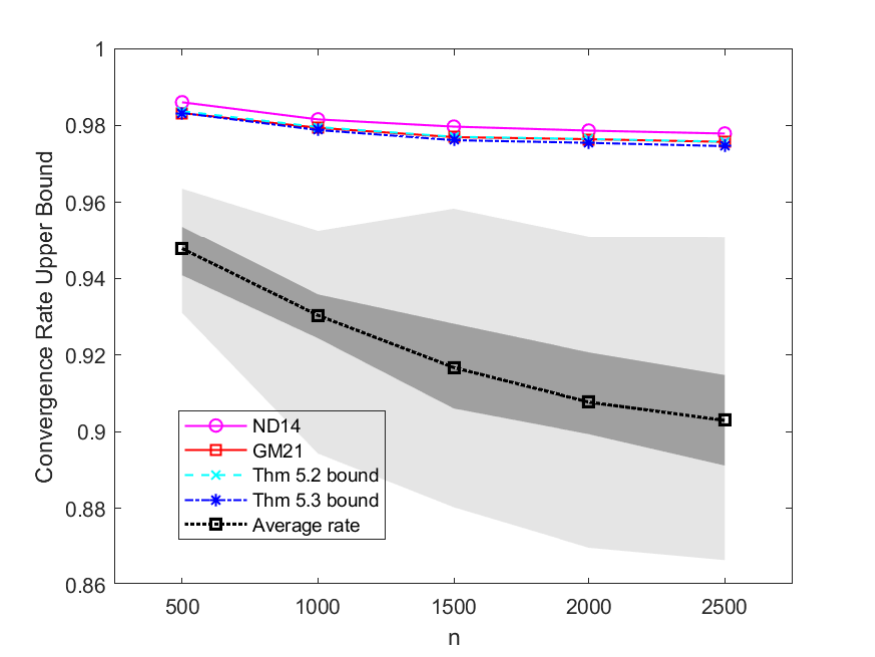}}
	\subfigure[]{\includegraphics[width=0.32\linewidth]{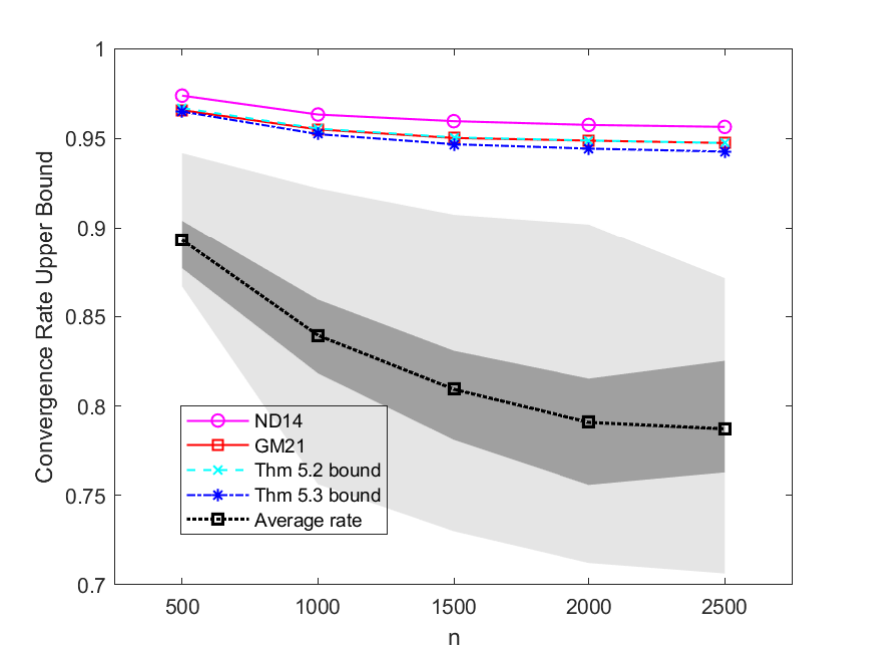}}
	\subfigure[]{\includegraphics[width=0.32\linewidth]{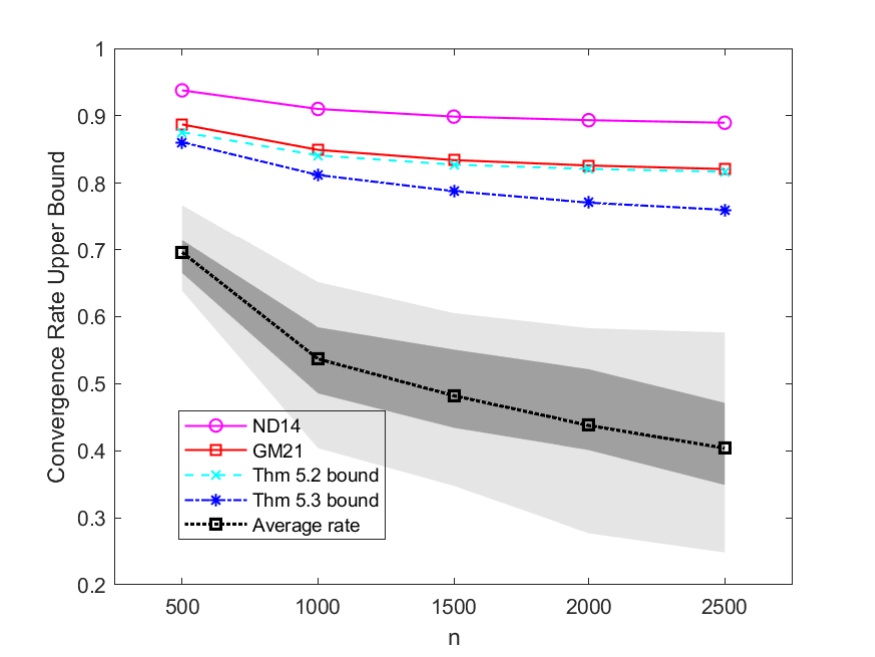}}
	\caption{Theoretical rate bounds and empirical convergence rate for different block sizes when $m=100$, tested on a two-scale matrix constructed from a Gaussian random matrix. Subfigures (a), (b) and (c) correspond to block sizes $q=10$, $q=20$ and $q=50$.}
	\label{fig:label1}
\end{figure}

For two-scale modification of randomly generated Gaussian matrix, FIG. \ref{fig:label1} show the comparison between the theoretical rate bounds and the empirical convergence rate under various block sizes $q$. We observe that the new theoretical rate bounds are consistently sharper than ND14. Furthermore, the  non-worst-case bound in Theorem \ref{convRBSK_sharper} is superior to GM21, and the gap becomes more significant as the block size $q$ and the ratio $n/m$ increase, while the worst-case result in Theorem \ref{convRBSK_relaxed_pave} yields a comparable bound to GM21.

\begin{figure}[H]
	\centering
	\subfigure[]{\includegraphics[width=0.32\linewidth]{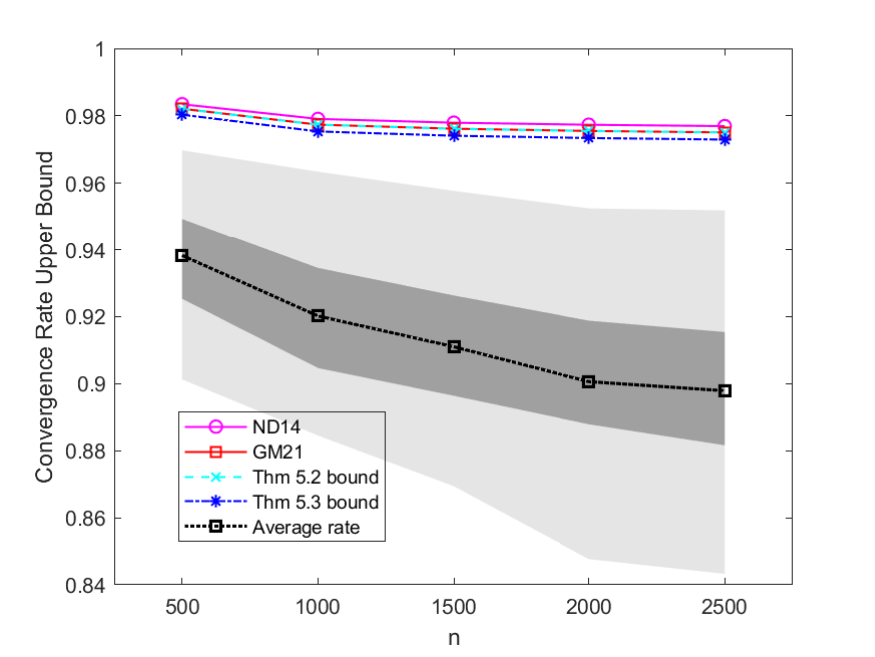}}
	\subfigure[]{\includegraphics[width=0.32\linewidth]{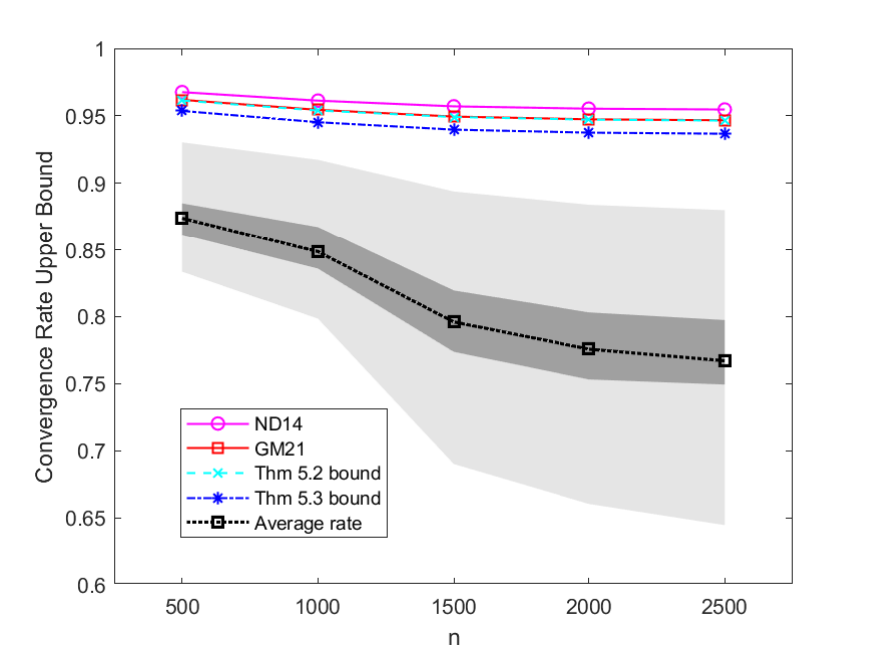}}
	\subfigure[]{\includegraphics[width=0.32\linewidth]{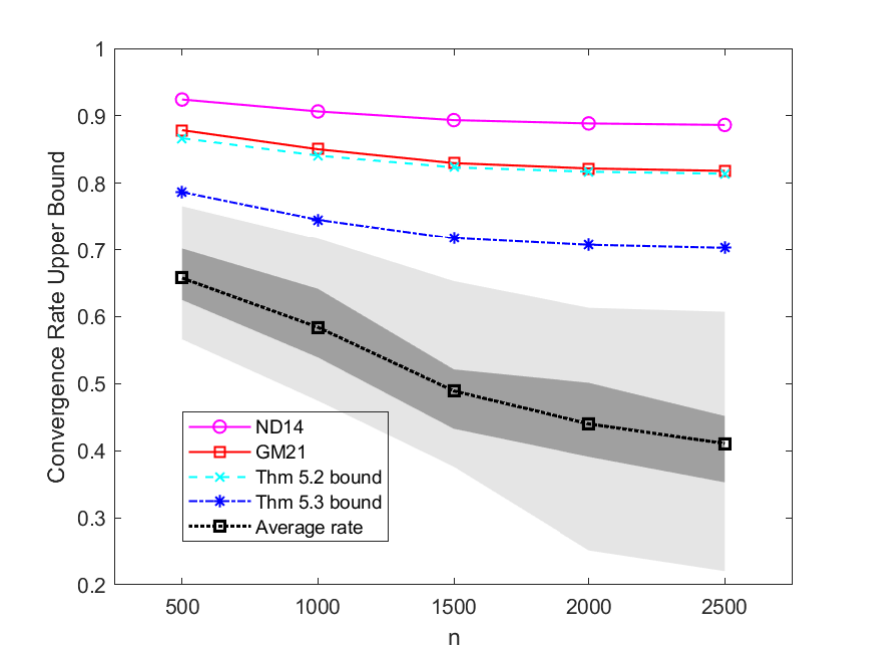}}
	\caption{Theoretical rate bounds and empirical convergence rate for different block sizes, tested on a two-scale matrix constructed from $n$ columns of \texttt{bibd\_81\_2}. Subfigures (a), (b) and (c) correspond to block sizes $q=10$, $q=20$ and $q=50$.}
	\label{fig:label3}
\end{figure}

\begin{figure}[H]
	\centering
	\subfigure[]{\includegraphics[width=0.32\linewidth]{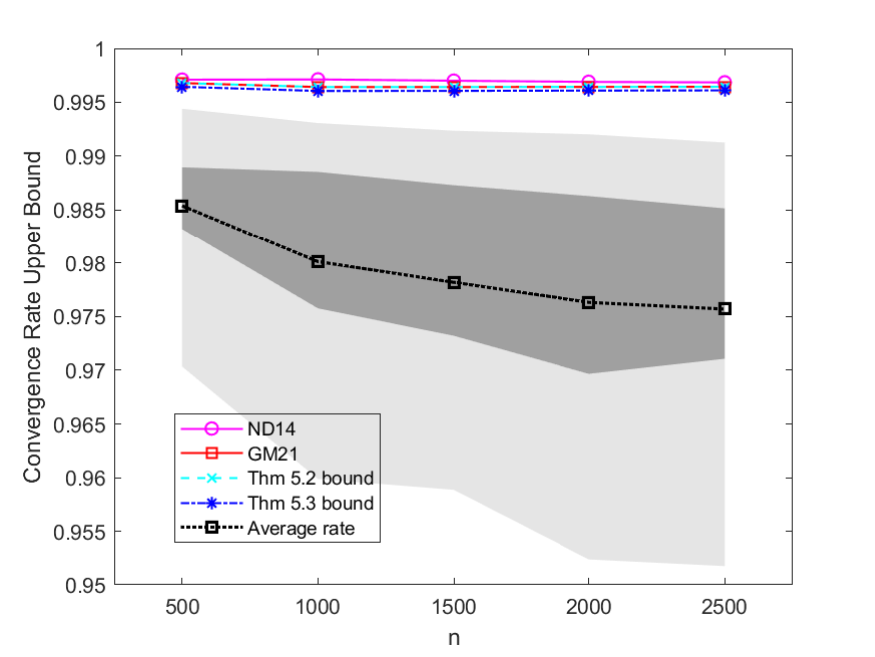}}
	\subfigure[]{\includegraphics[width=0.32\linewidth]{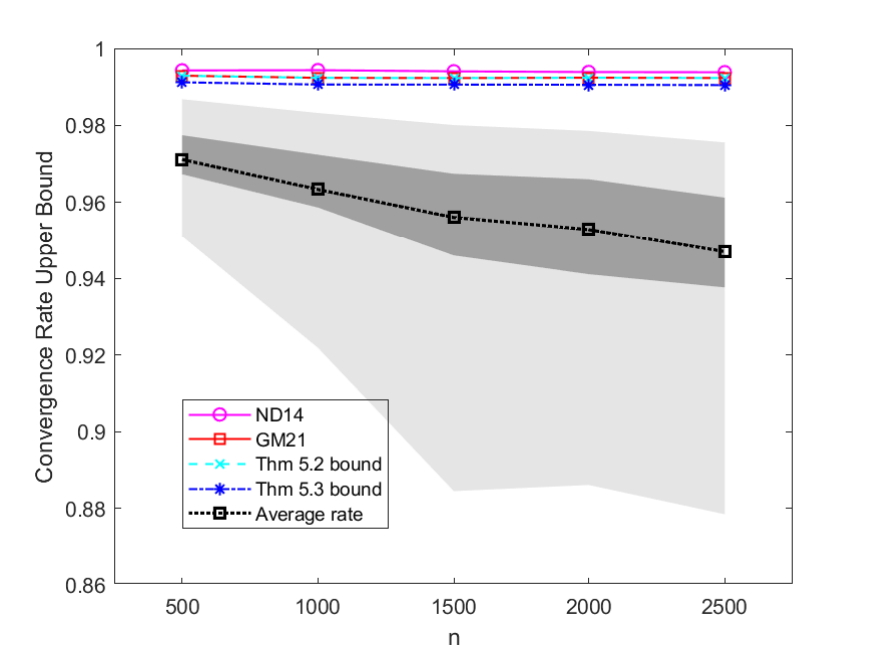}}
	\subfigure[]{\includegraphics[width=0.32\linewidth]{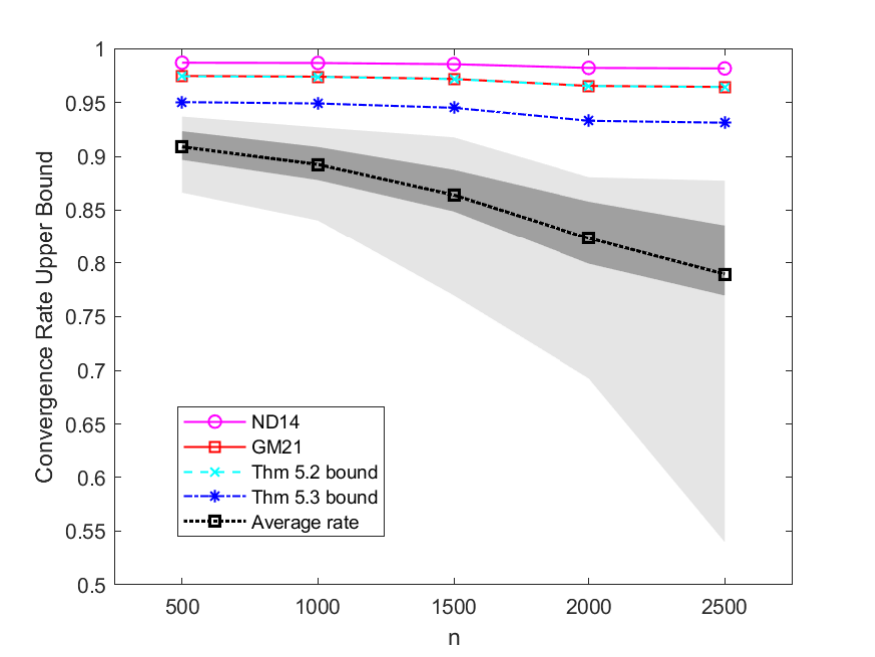}}
	\caption{Theoretical rate bounds and empirical convergence rate for different block sizes, tested on a two-scale matrix constructed from $n$ columns of \texttt{ch6-6-b5}. Subfigures (a), (b) and (c) correspond to block sizes $q=10$, $q=20$ and $q=50$.}
	\label{fig:label4}
\end{figure}

\begin{figure}[H]
	\centering
	\subfigure[]{\includegraphics[width=0.32\linewidth]{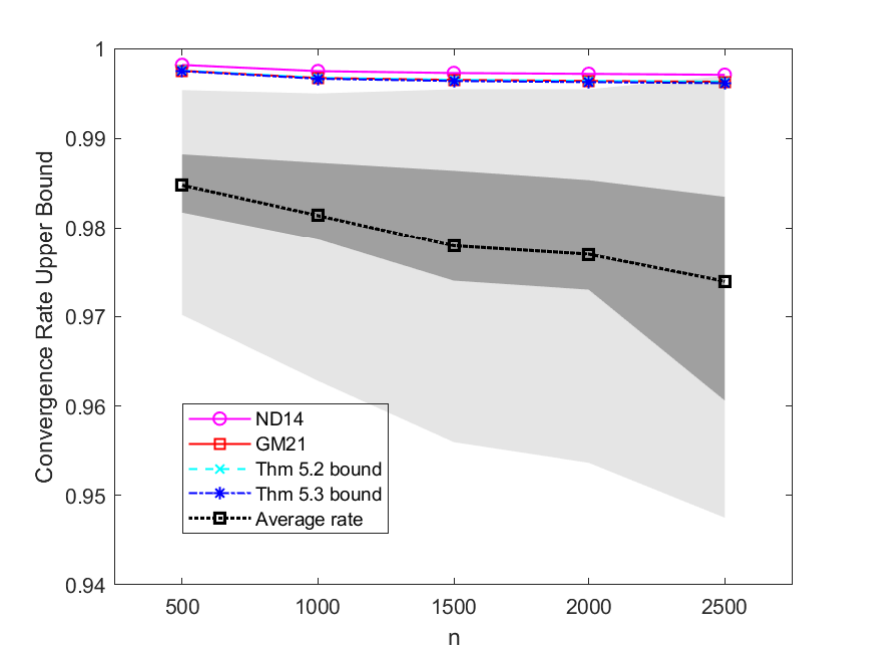}}
	\subfigure[]{\includegraphics[width=0.32\linewidth]{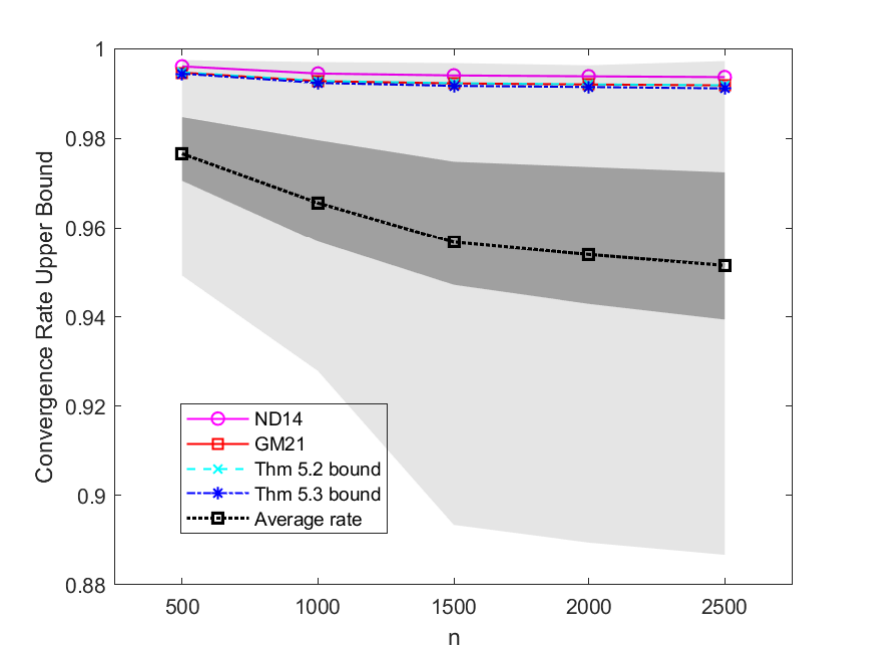}}
	\subfigure[]{\includegraphics[width=0.32\linewidth]{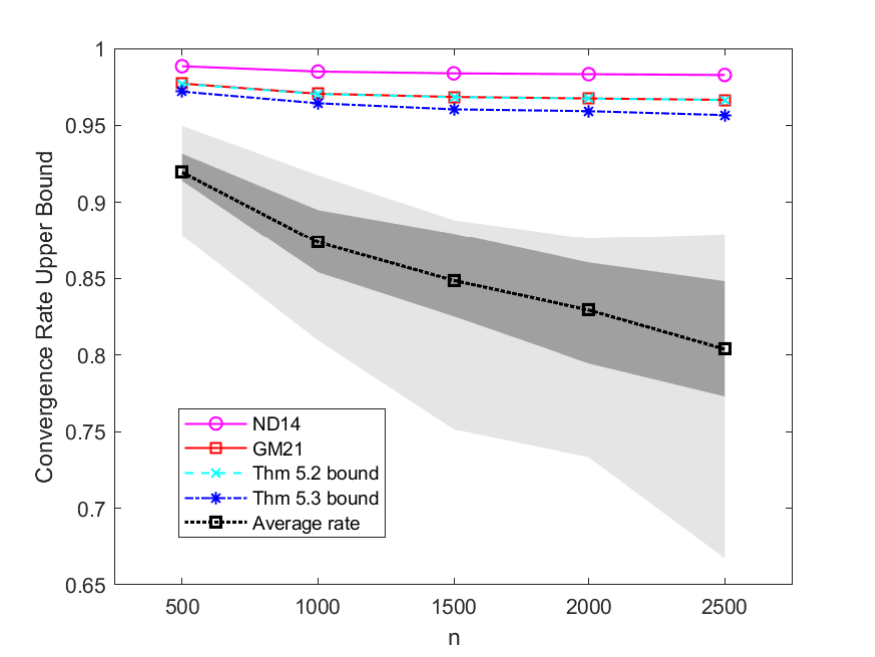}}
	\caption{Theoretical rate bounds and empirical convergence rate for different block sizes, tested on a two-scale matrix constructed from $n$ columns of \texttt{n4c5-b7}. Subfigures (a), (b) and (c) correspond to block sizes $q=10$, $q=20$ and $q=50$.}
	\label{fig:label5}
\end{figure}

Now, we consider real-world data matrices  \cite{davis2011university}. Three sparse matrices are tested: \texttt{bibd\_81\_2}, \texttt{ch6-6-b5}, and \texttt{n4c5-b7}. FIGs. \ref{fig:label3}-\ref{fig:label5} show the comparison between theoretical rate bounds and empirical convergence rates. For each matrix, a sub-data block is extracted and modified to introduce a two-scale structure. The results show that the bound in Theorem \ref{convRBSK_sharper} consistently outperforms both ND14 and GM21, and its advantage becomes more pronounced as the block size $q$ and the ratio $n/m$ increase. The bound in Theorem \ref{convRBSK_relaxed_pave} is generally comparable to GM21.

\subsection{Ill-conditioned matrices}
In this part, we construct ill-conditioned matrices by modifying a row block of the data matrix $A$ so that some singular values are significantly smaller than the larger ones. The construction proceeds as follows:
\begin{itemize}[itemsep=2pt, topsep=0pt, parsep=0pt, partopsep=0pt]
	\item compute the smallest positive singular value $\sigma_{\min}(A)$ of the data matrix $A$;
	\item generate Haar-distributed random orthogonal matrices $U \in \mathbb{R}^{q \times q}$ and $V \in \mathbb{R}^{n \times n}$;
	\item construct a diagonal matrix $\Sigma \in \mathbb{R}^{q \times n}$ with decreasing nonzero diagonal entries defined by $\sigma_i = \beta \sigma_{\min}(A) - (i-1)\ell$, where the scaling factor is $\beta=0.2$ and the decrement is $\ell=0.01$;
	\item form $B = U \Sigma V^\top$, and replace a row block of $A$ with $B$.
\end{itemize}

\begin{figure}[htbp]
	\centering
	\subfigure[]{\includegraphics[width=0.32\linewidth]{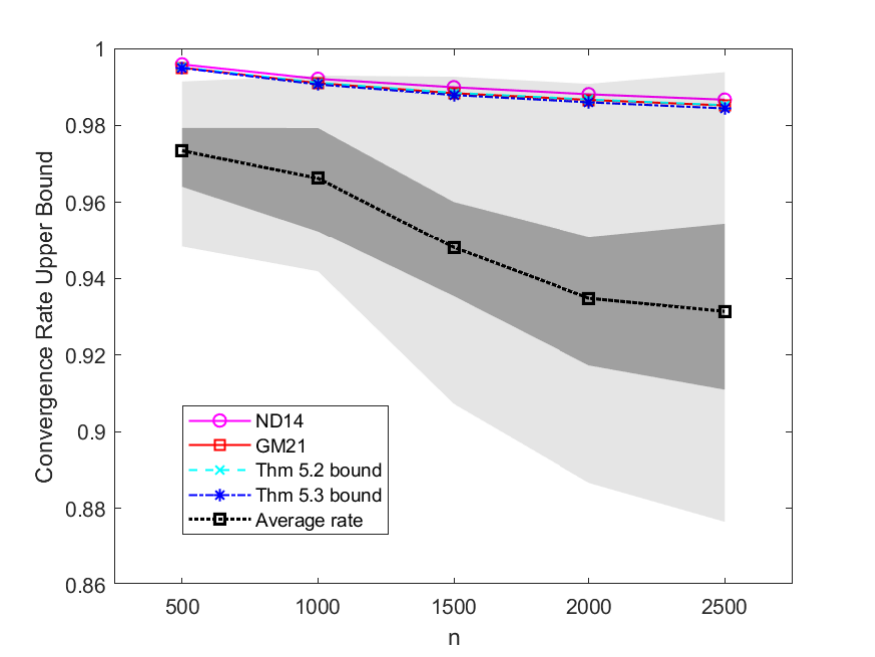}}
	\subfigure[]{\includegraphics[width=0.32\linewidth]{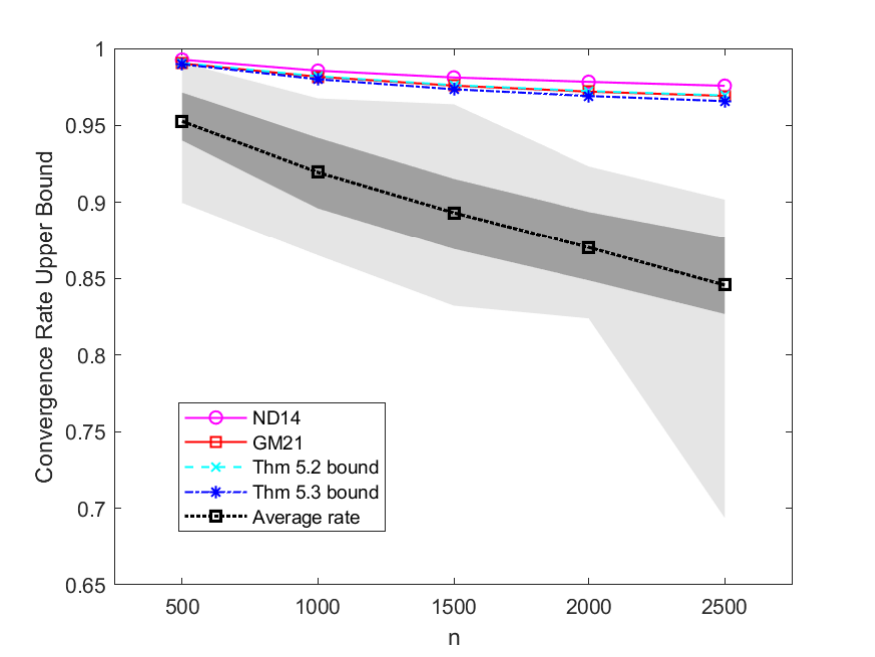}}
	\subfigure[]{\includegraphics[width=0.32\linewidth]{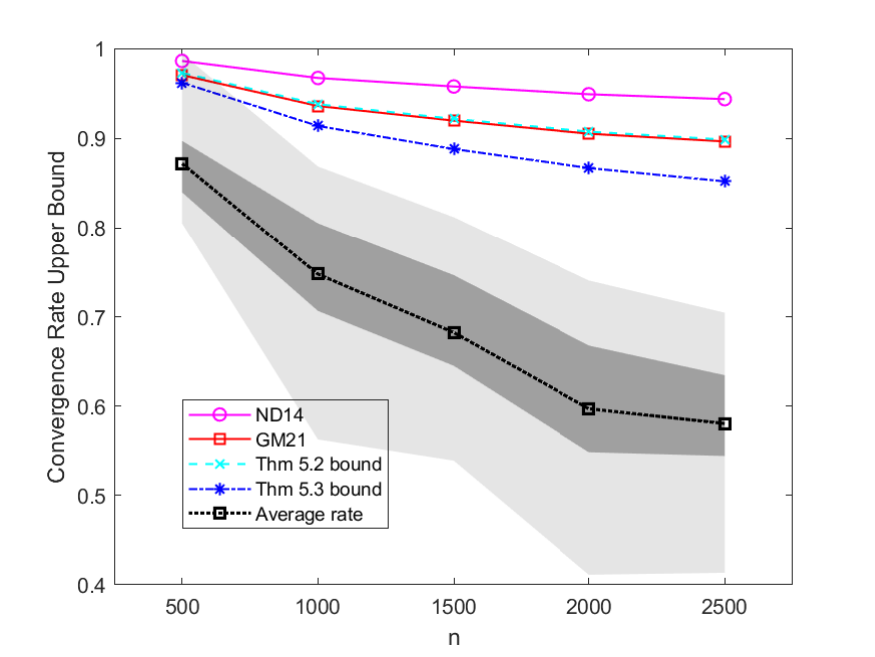}}
	\caption{Theoretical rate bounds and empirical convergence rate for different block sizes when $m=100$, tested on an ill-conditioned matrix constructed from a Gaussian random matrix. Subfigures (a), (b) and (c) correspond to block sizes $q=10$, $q=20$ and $q=50$.}
	\label{fig:label6}
\end{figure}

FIG. \ref{fig:label6} shows the comparison of theoretical and empirical convergence rates when the Gaussian matrix $A$ is randomly generated and modified in this way. As shown in FIG. \ref{fig:label6}, the new bounds consistently improve upon ND14. In particular, Theorem \ref{convRBSK_sharper} provides a significantly sharper non-worst-case bound than GM21, while the worst-case bound in Theorem \ref{convRBSK_relaxed_pave} is generally comparable to that of GM21.

\begin{figure}[H]
	\centering
	\subfigure[]{\includegraphics[width=0.32\linewidth]{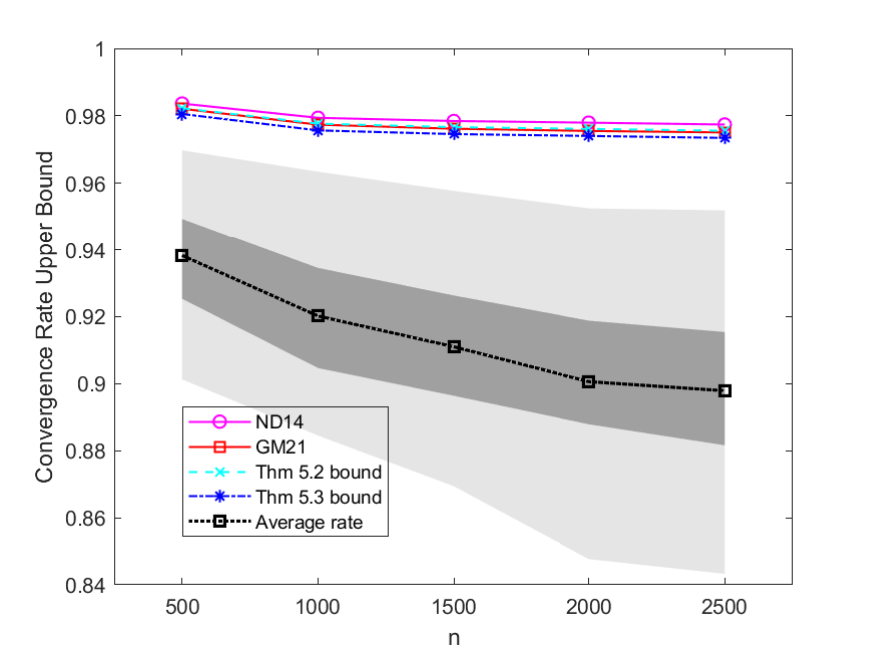}}
	\subfigure[]{\includegraphics[width=0.32\linewidth]{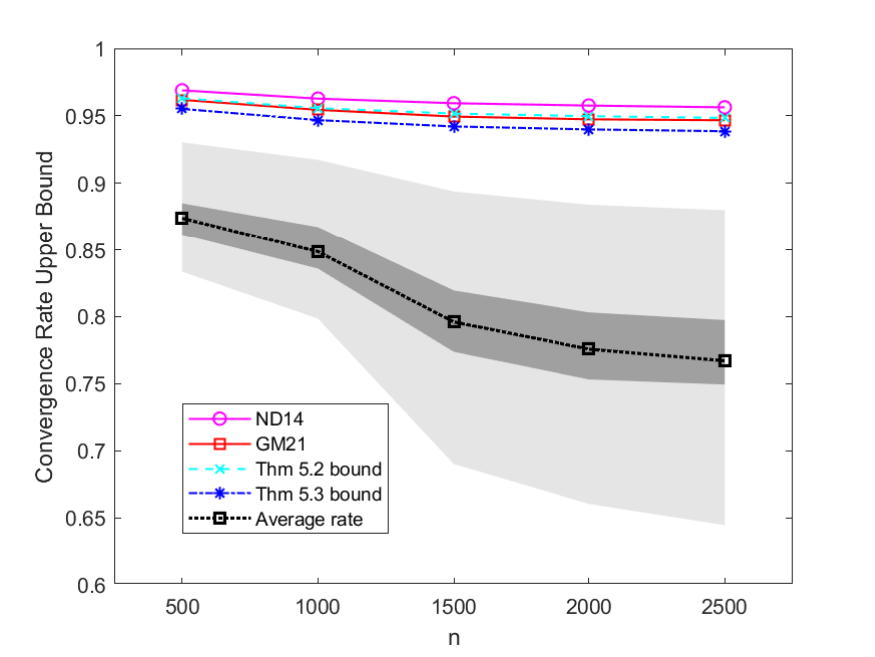}}
	\subfigure[]{\includegraphics[width=0.32\linewidth]{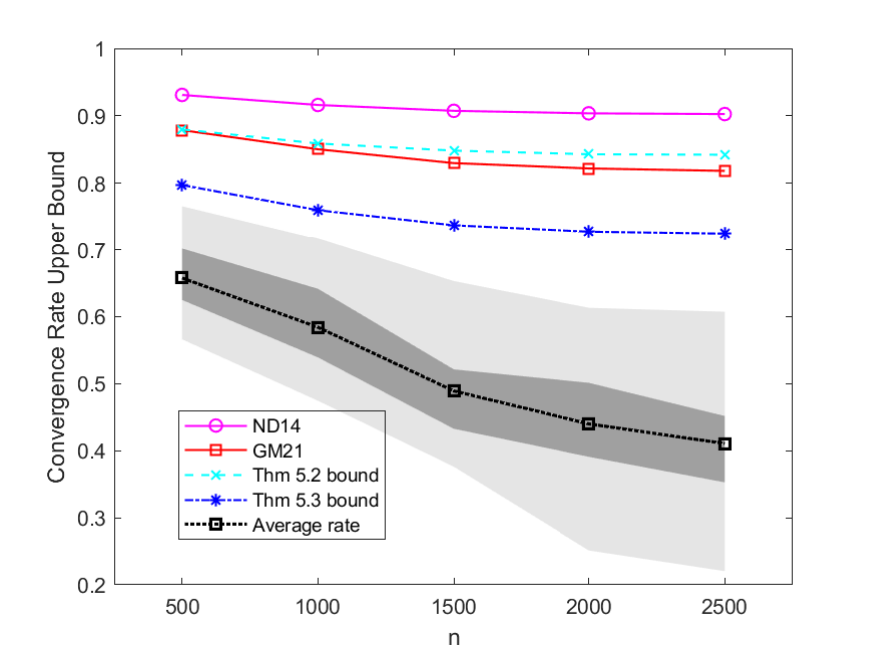}}
	\caption{Theoretical rate bounds and empirical convergence rate for different block sizes, tested on an ill-conditioned matrix constructed from $n$ columns of \texttt{bibd\_81\_2}. Subfigures (a), (b) and (c) correspond to block sizes $q=10$, $q=20$ and $q=50$.}
	\label{fig:label8}
\end{figure}

\begin{figure}[H]
	\centering
	\subfigure[]{\includegraphics[width=0.32\linewidth]{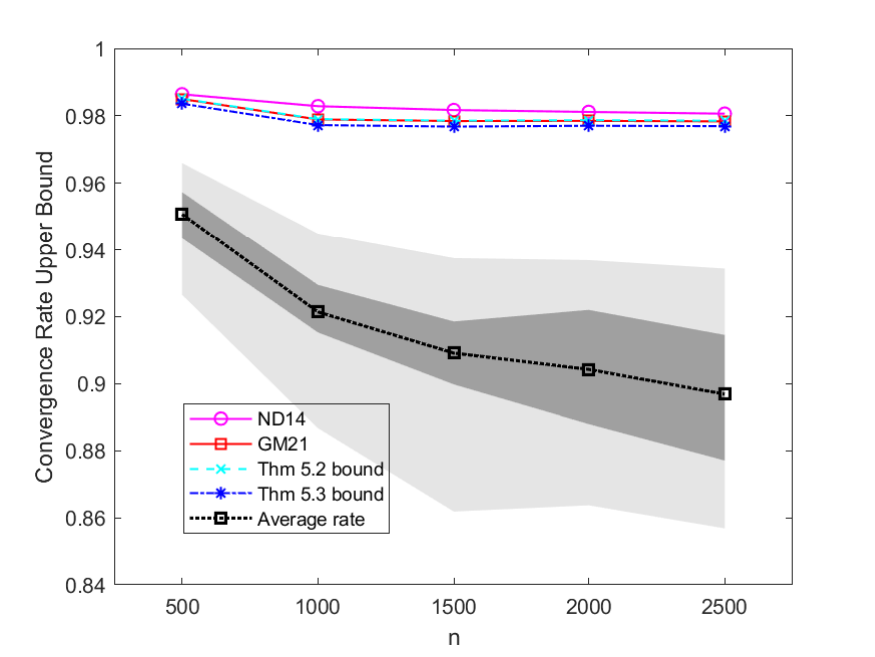}}
	\subfigure[]{\includegraphics[width=0.32\linewidth]{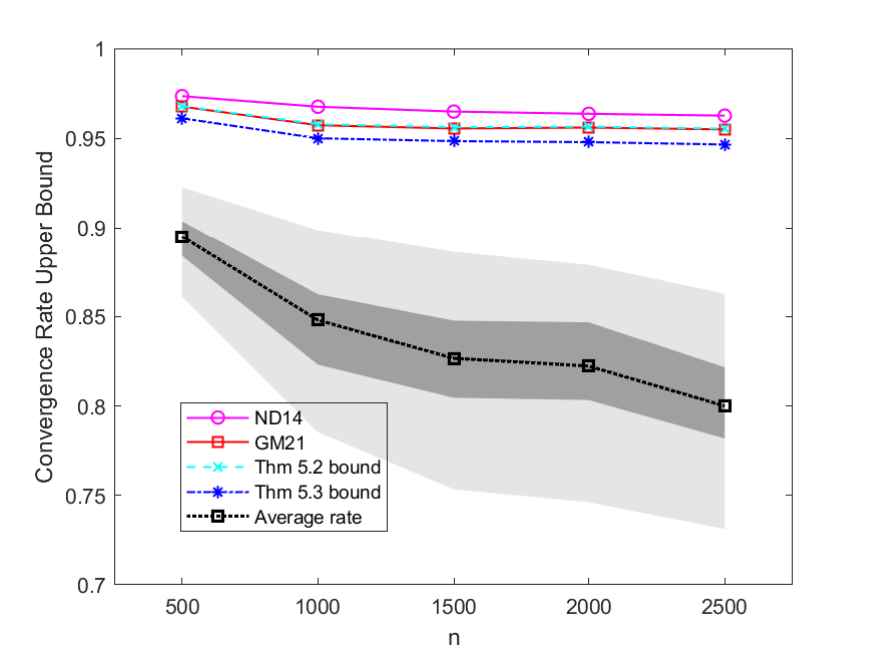}}
	\subfigure[]{\includegraphics[width=0.32\linewidth]{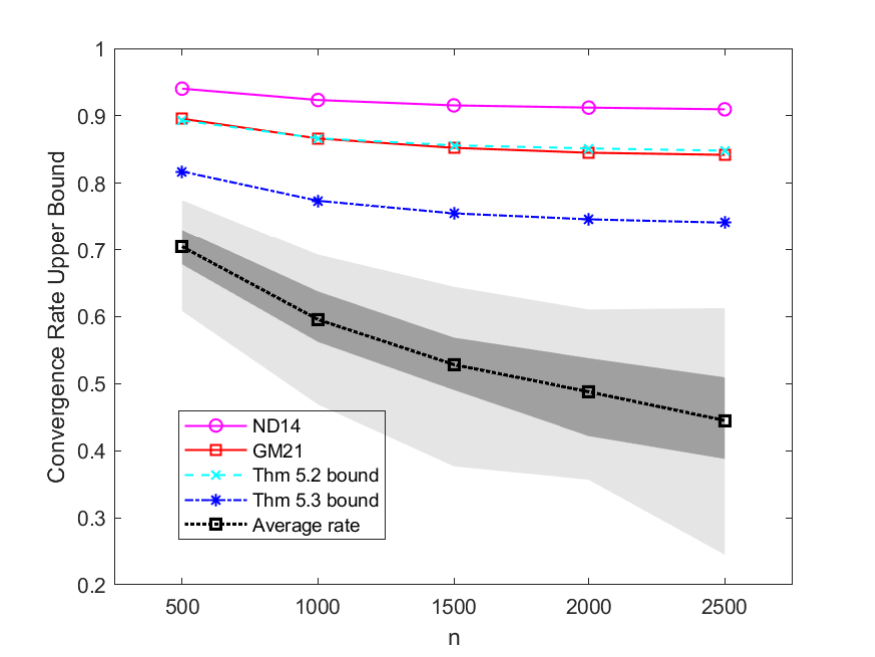}}
	\caption{Theoretical rate bounds and empirical convergence rate for different block sizes, tested on an ill-conditioned matrix constructed from $n$ columns of \texttt{ch6-6-b5}. Subfigures (a), (b) and (c) correspond to block sizes $q=10$, $q=20$ and $q=50$.}
	\label{fig:label9}
\end{figure}

\begin{figure}[H]
	\centering
	\subfigure[]{\includegraphics[width=0.32\linewidth]{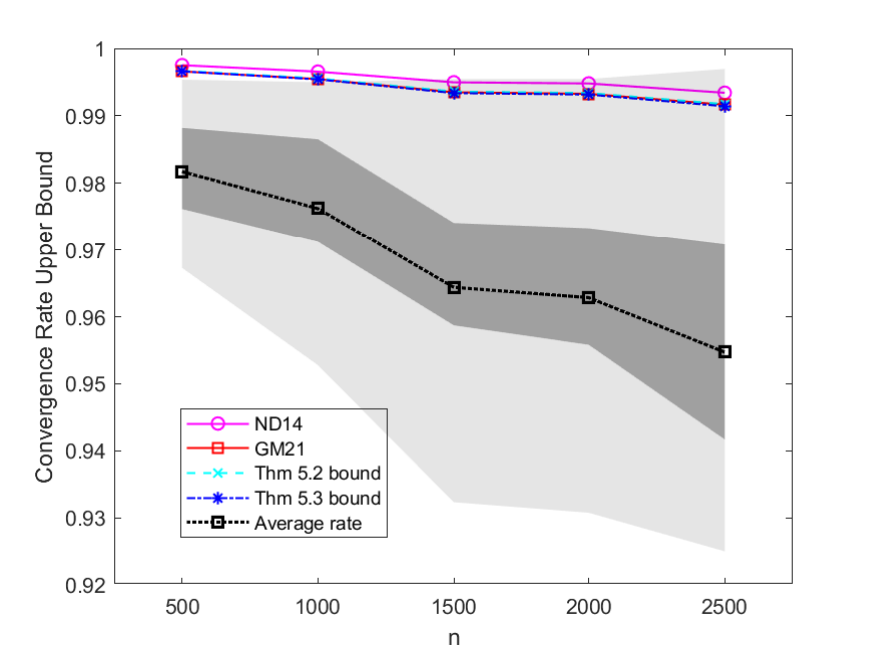}}
	\subfigure[]{\includegraphics[width=0.32\linewidth]{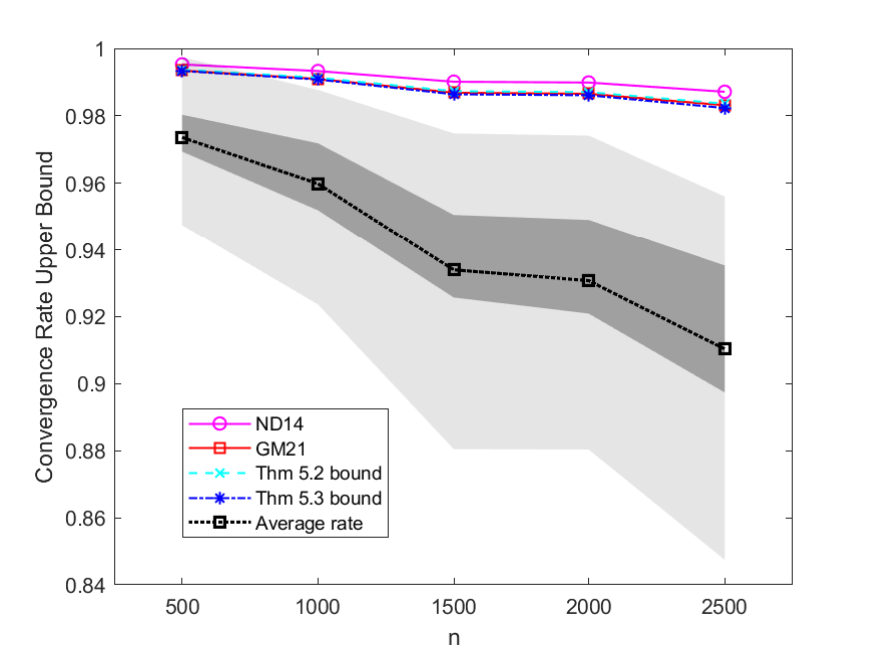}}
	\subfigure[]{\includegraphics[width=0.32\linewidth]{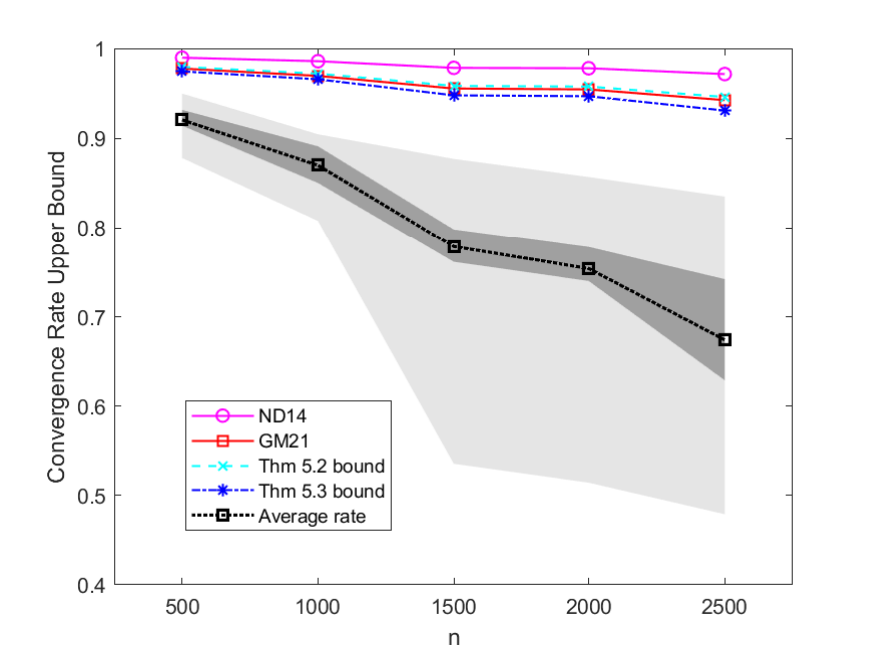}}
	\caption{Theoretical rate bounds and empirical convergence rate for different block sizes, tested on an ill-conditioned matrix constructed from $n$ columns of \texttt{n4c5-b7}. Subfigures (a), (b) and (c) correspond to block sizes $q=10$, $q=20$ and $q=50$.}
	\label{fig:label10}
\end{figure}

We further consider test problems constructed from real-world data matrices \cite{davis2011university}. FIGs. \ref{fig:label8}-\ref{fig:label10} present results for the modified matrices
\texttt{bibd\_81\_2}, \texttt{ch6-6-b5}, and \texttt{n4c5-b7}. Across all three matrices, the bound of Theorem \ref{convRBSK_sharper} consistently outperforms GM21 and ND14, with its advantage becoming more pronounced as the block size $q$ and the ratio $n/m$ increase, while the bound of Theorem \ref{convRBSK_relaxed_pave} is generally comparable to GM21.

\section{Conclusions}
\label{sec:concl}

In this paper, we develop a unified analysis technique of randomized Kaczmarz-type methods through the RBSK framework. Using concentration inequalities, we obtain new tight expected linear convergence rate bounds. The introduced scaling operator $S$ makes the new bounds scale-invariant, thereby eliminating the dependence on the magnitude of the data matrix. In addition, we also clarify the connection between our new bounds and existing bounds for randomized non-extended block Kaczmarz methods, see Remarks \ref{remark:relaxOne} and \ref{remark:relaxTwo}.

Numerical experiments on synthetic multi-scale and ill-conditioned instances, as well as on sparse matrices from SuiteSparse, show that the proposed bounds are tight to the empirical convergence rates. In most cases, the refined new bounds are sharper than ND14 and GM21, and this advantage becomes more evident as the block size $q$ increases.

The present work addresses randomized block Kaczmarz with static sampling rules, but several natural extensions remain. Development and analysis of learning-guided batch-sampling distributions, as well as other approaches for distribution optimization, are valuable topics for future research. The new analytical technique may be further developed to handle iterate-dependent sampling strategies, such as greedy and adaptive variants. In addition, similar analysis techniques could be applied to study the theoretical convergence rates of randomized batch-sampling coordinate descent methods, thereby broadening the applicability of the new framework.

\section*{Acknowledgments}
The corresponding author is supported by the National Natural Science Foundation of China (Grant Nos. 12071215, 11101213).

\appendix
\section{Other Lemmas}

\label{otherLemmas}

\begin{lemma}
  \label{lem2}
  Let $A\in\reals^{m\times n}$ be a nonzero matrix. If $u\in\range(A)$, it holds that
  \begin{eqnarray*}
    \|A^{\dag}u\|_2^2 & \geq & \frac{1}{\|A\|_2^2}\|u\|_2^2.
  \end{eqnarray*}
\end{lemma}
{\em Proof.}
Let $\rank(A)=s>0$. Then $A$ admits the following reduced singular value decomposition
\begin{eqnarray*}
  A &=& V_s\Sigma_sU_s^{\T},
\end{eqnarray*}
where $V_s\in\reals^{m\times s}$ and $U_s\in\reals^{n\times s}$ have
orthonormal columns, and $\Sigma_s=(\sigma_1,\ldots,\sigma_s)\in\reals^{s\times s}$
is a diagonal matrix with decreasing positive singular values of $A$. Obviously, there is a fact $\range(A)=\range(V_s)$, which leads to
\begin{eqnarray*}
  u &=& V_su_s \mbox{ for some } u_s\in\reals^s, \mbox{ and } \|u\|_2=\|u_s\|_2.
\end{eqnarray*}
Since the Moore-Penrose pseudoinverse of $A$ reads
\begin{eqnarray*}
  A^{\dag} &=& U_s\Sigma_s^{-1}V_s^{\T},
\end{eqnarray*}
then it follows that
\begin{eqnarray*}
  A^{\dag}u &=& U_s\Sigma_s^{-1}u_s.
\end{eqnarray*}
Therefore, it holds that
\begin{eqnarray*}
  \|A^{\dag}u\|_2^2 &=& \|U_s\Sigma_s^{-1}u_s\|_2^2 \\
                    &=& \|\Sigma_s^{-1}u_s\|_2^2 \\
                    &\geq& \frac{1}{\sigma_1^2}\|u_s\|_2^2.
\end{eqnarray*}
Together with the facts $\|u\|_2=\|u_s\|_2$ and $\sigma_1=\|A\|_2$, one may obtain the result of Lemma \ref{lem2}. $\hfill\square$

\begin{lemma}
  \label{lem3}
  Let $A\in\reals^{m\times n}$ be a nonzero matrix, $y\in\range(A)$ and $x\in\reals^n$ be nonzero vectors.
  Let $S\in\reals^{m\times m}$ be a nonsingular matrix, then
  the linear systems $Ax=y$ and $SAx=Sy$ have the same
  least-norm solution, i.e.,
  \begin{eqnarray*}
    A^{\dag}y &=& (SA)^{\dag}Sy.
  \end{eqnarray*}
\end{lemma}
{\em Proof.}
Due to the fact $y\in\range(A)$, the solution space of
the linear system $Ax=y$ is nonempty. Since the matrix
$S$ is nonsingular, the solution space of $SAx=Sy$ is the
same as that of $Ax=y$. Therefore, the least-norm solution
of $Ax=y$ and $SAx=Sy$ are the same. $\hfill\square$

\begin{lemma}
  \label{lem4}
  Let $X$ and $Y$ be random variables with $\E(|X|+|Y|)<\infty$, then
  \begin{eqnarray*}
    \E(XY) &=& \E\left[X \E(Y|X)\right].
  \end{eqnarray*}
\end{lemma}
{\em Proof.} This result is a direct consequence of the law of total expectation \cite{ChungKL2001B}.

\begin{lemma}[\cite{BaiWu18SIAM,DuSiSun20SISC}]
  \label{lem1}
  Let $A\in\reals^{m\times n}$ be a nonzero matrix. If $u\in\range(A^{\T})$, it holds that
  \begin{eqnarray*}
    \|Au\|_2^2 & \geq & \sigma_{\min}^2(A)\|u\|_2^2.
  \end{eqnarray*}
\end{lemma}

\begin{lemma}
	\label{lem5}
	Let $\{x_i\}_{i=1}^n$ and $\{y_i\}_{i=1}^n$ be two sample sequences of size $n$, with sample means $\bar{x}$ and $\bar{y}$ respectively. If the sample covariance satisfies
	\begin{eqnarray*}
		\cov\left(x, y\right)=\frac{1}{n-1} \sum_{i=1}^n\left(x_i-\bar{x}\right)\left(y_i-\bar{y}\right) \geq 0,
	\end{eqnarray*}
	then
	\begin{eqnarray*}
		\frac{1}{n}\sum_{i=1}^n x_i y_i \geq \bar{x}\bar{y}=\left(\frac{1}{n}\sum_{i=1}^n x_i\right) \left(\frac{1}{n}\sum_{i=1}^n y_i\right).
	\end{eqnarray*}
\end{lemma}

\section{Examples for Batch-Sampling Rules}

\label{examples_samplingRules}

\begin{example}\label{exampRowPaving}
  \textsc{Row paving batch-sampling}. When the range of the batch-sampling $\tau$ with a prescribed joint distribution $\bbb{P}$
  refers to a row paving of the data matrix $A\in\reals^{m\times n}$, the batch-sampling $\tau\sim\bbb{P}$
  reduces to the case studied in {\rm \cite{NeedellTropp2014LAA}}, see also the \textsc{partition sampling} case
  studied in {\rm \cite{Necoara19SIMAX}}. In particular,
  a ($L,\beta_{\txt{\tiny Low}},\beta_{\txt{\tiny Up}}$) row paving of the data matrix $A$ is introduced
  in {\rm \cite{NeedellTropp2014LAA}} by defining a partition $\mathcal{T}=\{T_1,\ldots,T_{L}\}$ of the row indices $[m]$ that satisfies
 \begin{eqnarray}
  \label{def:rowPaving}
    \beta_{\txt{\tiny Low}}\leq\lambda_{\min}(A_{\tau}A_{\tau}^{\T})\ \mbox{ and }\
    \lambda_{\max}(A_{\tau}A_{\tau}^{\T})\leq\beta_{\txt{\tiny Up}}\ \mbox{ for each }\
    \tau\in \mathcal{T}.
  \end{eqnarray}
  The RBK method proposed in {\rm \cite{NeedellTropp2014LAA}} is designed by selecting a batch-sampling $\tau$ at random at each iteration such that the effective sets
  $\tau_{\txttiny{E}}\in\range(\tau)=\mathcal{T}$ are
  uniformly (equal probability $\frac{1}{L}$) or non-uniformly (unequal probability) sampled.
  Let $q=\max_{1\leq i\leq L}|T_i|$, the joint distribution $\bbb{P}$ of the batch-sampling $\tau$
  can be represented by a $q$th-order tensor with constant length (i.e., $m$) of each dimension.
  For instance, let's consider a data matrix $A\in\reals^{3\times n}$, which admits a row paving
  ($2,\beta_{\mbox{\tiny Low}},\beta_{\mbox{\tiny Up}}$) with partition $\mathcal{T}=\{T_1=\{1\},T_2=\{2,3\}\}$
  and $q=\max \{|T_1|=1,|T_2|=2\}=2$. In addition,
  the RBK method selects a batch-sampling $\tau=(\tau_1,\tau_2)$ at random at each iteration
  such that the effective sets $T_1=\{1\}$ and $T_2=\{2,3\}$ are uniformly (equal probability $\frac{1}{2}$) sampled,
  i.e.,
  \begin{eqnarray*}
    \pr(\{1\})	&=& \pr(\tau=(1,1)) = \bbb{p}_{11} = \frac{1}{2}, \\
    \pr(\{2,3\})  &=& \pr(\tau=(2,3))+\pr(\tau=(3,2)) = \bbb{p}_{23}+\bbb{p}_{32} = \frac{1}{2}.
  \end{eqnarray*}
  Then, the batch-sampling $\tau$ satisfies the joint distribution $\bbb{P}$ listed in
  Table \ref{tab:jointDistribRowPavingSampling-3-by-3}, which is a 2nd-order tensor with constant length 3 of each dimension,
  i.e., a 3-by-3 square matrix.
  The marginal distributions of $\tau_1$ and $\tau_2$ are listed at the last column and row of
  Table \ref{tab:jointDistribRowPavingSampling-3-by-3}, respectively.
\renewcommand\arraystretch{1.5}
\begin{table}[htbp]
\setlength{\abovecaptionskip}{0pt}
\setlength{\belowcaptionskip}{10pt} \centering{
\caption{\label{tab:jointDistribRowPavingSampling-3-by-3}
The joint distribution $\bbb{P}$ of row paving batch-sampling $\tau=(\tau_1,\tau_2)$,
and the marginal distributions of $\tau_1$ and $\tau_2$: $A\in\reals^{3\times n}$, $\mathcal{T}=\{T_1=\{1\},T_2=\{2,3\}\}$.}
\begin{tabular}{|c|ccc|c|}\hline
\diagbox{$\tau_1$}{$\tau_2$} & 1 & 2 & 3 & $\tau_1=i$ \\ \hline
                           1 & $\frac{1}{2}$ & 0 & 0 & $\frac{1}{2}$ \\
                           2 & 0 & 0 & $\bbb{p}_{23}$ & $\bbb{p}_{23}$ \\
                           3 & 0 & $\frac{1}{2}-\bbb{p}_{23}$ & 0 & $\frac{1}{2}-\bbb{p}_{23}$ \\ \hline
                  $\tau_2=j$ & $\frac{1}{2}$ & $\frac{1}{2}-\bbb{p}_{23}$ & $\bbb{p}_{23}$ & 1 \\
\hline
\end{tabular}}
\end{table}

According to Definition \ref{def:batchSampling}, the diagonal matrices $P_1$ and $P_2$ for $\tau_1$ and
  $\tau_2$ are of the forms
  \begin{eqnarray*}
    P_1 = \begin{bmatrix}
              \frac{1}{2} &   &   \\
                & \bbb{p}_{23} &   \\
                &   & \frac{1}{2}-\bbb{p}_{23} \\
            \end{bmatrix},\
    P_2 = \begin{bmatrix}
              \frac{1}{2} &   &   \\
                & \frac{1}{2}-\bbb{p}_{23} &   \\
                &   & \bbb{p}_{23} \\
            \end{bmatrix}.
  \end{eqnarray*}
  Obviously, the property $P_1+P_2\succ 0$ gets satisfied.

\end{example}

\begin{example}\label{exampUniform}
  \textsc{Uniform batch-sampling}. When the range of the batch-sampling $\tau$ with a prescribed
  joint distribution $\bbb{P}$ defines $\range(\tau)=\mathcal{F}\subseteq 2^{[m]}$ (power set of $[m]$),
  and each effective set $\tau_{\txttiny{E}}\in\mathcal{F}$ includes $q$ unique row indices,
  i.e., $n_{\txttiny{E}}=q$ for all $\tau_{\txttiny{E}}\in\mathcal{F}$,
  there are $\tbinom{m}{q}$ possible values of the effective sets $\tau_{\txttiny{E}}\in\mathcal{F}$,
  thus, the batch-sampling $\tau$ reduces to one of the cases handled by {\rm \cite{Necoara19SIMAX}}.
  One case of
  the randomized average block Kaczmarz (RaBK) method proposed in {\rm \cite{Necoara19SIMAX}}
  is designed by selecting a batch-sampling $\tau$ at random at each iteration such that
  the effective sets $\tau_{\txttiny{E}}\in\mathcal{F}$ are uniformly (equal probability $1/\tbinom{m}{q}$)
  or non-uniformly (unequal probability) sampled,
  which can be considered as a pseudo-inverse free randomized variant of the block Kaczmarz method.
  For instance, let's consider a data matrix $A\in\reals^{3\times n}$, and $q=2$, then there are
  3 possible values of the effective sets $\tau_{\txttiny{E}}\in\mathcal{F}=\{\{1,2\},\{1,3\},\{2,3\}\}$.
  In addition,
  the RaBK method selects a batch-sampling $\tau=(\tau_1,\tau_2)$ at random at
  each iteration such that the effective sets $\{1,2\}$, $\{1,3\}$ and $\{2,3\}$
  are uniformly (equal probability $\frac{1}{3}$) sampled, i.e.,
  \begin{eqnarray*}
    \pr(\{1,2\}) &=& \pr(\tau=(1,2))+\pr(\tau=(2,1)) = \bbb{p}_{12}+\bbb{p}_{21} = \frac{1}{3}, \\
    \pr(\{1,3\}) &=& \pr(\tau=(1,3))+\pr(\tau=(3,1)) = \bbb{p}_{13}+\bbb{p}_{31} = \frac{1}{3}, \\
    \pr(\{2,3\}) &=& \pr(\tau=(2,3))+\pr(\tau=(3,2)) = \bbb{p}_{23}+\bbb{p}_{32} = \frac{1}{3}.
  \end{eqnarray*}
  Then, the batch-sampling $\tau$ satisfies the joint distribution $\bbb{P}$ listed in
  Table \ref{tab:jointDistribUniformSampling-3-by-3}, which is again a 2nd-order tensor with constant length 3 of each dimension.
  The marginal distributions of $\tau_1$ and $\tau_2$ are listed at the last column and row of
  Table \ref{tab:jointDistribUniformSampling-3-by-3}, respectively.
\renewcommand\arraystretch{1.5}
\begin{table}[htbp]
\setlength{\abovecaptionskip}{0pt}
\setlength{\belowcaptionskip}{10pt} \centering{
\caption{\label{tab:jointDistribUniformSampling-3-by-3}
The joint distribution $\bbb{P}$ of uniform batch-sampling $\tau=(\tau_1,\tau_2)$,
and the marginal distributions of $\tau_1$ and $\tau_2$: $A\in\reals^{3\times n}$, $\mathcal{F}=\{\{1,2\},\{1,3\},\{2,3\}\}$.}
\begin{tabular}{|c|ccc|c|}\hline
\diagbox{$\tau_1$}{$\tau_2$} & 1 & 2 & 3 & $\tau_1=i$ \\ \hline
                           1 & 0 & $\bbb{p}_{12}$ & $\bbb{p}_{13}$ & $\bbb{p}_{12}+\bbb{p}_{13}$ \\
                           2 & $\frac{1}{3}-\bbb{p}_{12}$ & 0 & $\bbb{p}_{23}$ & $\frac{1}{3}-\bbb{p}_{12}+\bbb{p}_{23}$ \\
                           3 & $\frac{1}{3}-\bbb{p}_{13}$ & $\frac{1}{3}-\bbb{p}_{23}$ & 0 & $\frac{2}{3}-\bbb{p}_{13}-\bbb{p}_{23}$ \\ \hline
                  $\tau_2=j$ & $\frac{2}{3}-\bbb{p}_{12}-\bbb{p}_{13}$ & $\frac{1}{3}+\bbb{p}_{12}-\bbb{p}_{23}$ & $\bbb{p}_{13}+\bbb{p}_{23}$ & 1 \\
\hline
\end{tabular}}
\end{table}
  According to Definition \ref{def:batchSampling}, the diagonal matrices $P_1$ and $P_2$ for $\tau_1$ and
  $\tau_2$ are of the forms
  \begin{eqnarray*}
    P_1 = \begin{bmatrix}
              \bbb{p}_{12}+\bbb{p}_{13} &   &   \\
                & \frac{1}{3}-\bbb{p}_{12}+\bbb{p}_{23} &   \\
                &   & \frac{2}{3}-\bbb{p}_{13}-\bbb{p}_{23} \\
            \end{bmatrix},\
    P_2 = \begin{bmatrix}
              \frac{2}{3}-\bbb{p}_{12}-\bbb{p}_{13} &   &   \\
                & \frac{1}{3}+\bbb{p}_{12}-\bbb{p}_{23} &   \\
                &   & \bbb{p}_{13}+\bbb{p}_{23} \\
            \end{bmatrix}.
  \end{eqnarray*}
  Obviously, the property $P_1+P_2\succ 0$ also gets satisfied.
\end{example}

\begin{example}\label{exampNonUnique}
  \textsc{Non-unique batch-sampling}. When the range of the batch-sampling $\tau$ with a prescribed
  joint distribution $\bbb{P}$ defines $\range(\tau)=\mathcal{G}\subseteq 2^{[m]}$,
  and each effective sets $\tau_{\txttiny{E}}\in\mathcal{G}$ includes at most $q$ unique row indices,
  i.e., $1\leq n_{\txttiny{E}}\leq q$ for all $\tau_{\txttiny{E}}\in\mathcal{G}$,
  there are $\sum_{t=1}^{q}\tbinom{m}{t}$ possible values of the effective sets $\tau_{\txttiny{E}}\in\mathcal{G}$,
  thus, the batch-sampling $\tau$ refers to a new stochastic sampling that differs from the row paving and uniform cases
  in Examples \ref{exampRowPaving} and \ref{exampUniform}.
  Then, a RBSK method can be obtained
  by selecting a batch-sampling $\tau$ at random at each iteration such that the effective sets $\tau_{\txttiny{E}}\in\mathcal{G}$
  are uniformly (equal probability $1/\sum_{t=1}^{q}\tbinom{m}{t}$) or non-uniformly (unequal probability) sampled.
  For instance, let's consider a data matrix $A\in\reals^{3\times n}$, and $q=2$, then there are
  6 possible values of the effective sets $\tau_{\txttiny{E}}\in\mathcal{G}=\{\{1\},\{2\},\{3\},\{1,2\},\{1,3\},\{2,3\}\}$.
  In addition, the RBSK method selects a batch-sampling $\tau=(\tau_1,\tau_2)$ at random at each iteration
  such that the effective sets $\{1\}$, $\{2\}$, $\{3\}$, $\{1,2\}$, $\{1,3\}$ and $\{2,3\}$
  are uniformly (equal probability $\frac{1}{6}$) sampled, i.e.,
  \begin{eqnarray*}
    \pr(\{1\}) &=& \pr(\tau=(1,1)) = \bbb{p}_{11} = \frac{1}{6}, \\
    \pr(\{2\}) &=& \pr(\tau=(2,2)) = \bbb{p}_{22} = \frac{1}{6}, \\
    \pr(\{3\}) &=& \pr(\tau=(3,3)) = \bbb{p}_{33} = \frac{1}{6}, \\
    \pr(\{1,2\}) &=& \pr(\tau=(1,2))+\pr(\tau=(2,1)) = \bbb{p}_{12}+\bbb{p}_{21} = \frac{1}{6}, \\
    \pr(\{1,3\}) &=& \pr(\tau=(1,3))+\pr(\tau=(3,1)) = \bbb{p}_{13}+\bbb{p}_{31} = \frac{1}{6}, \\
    \pr(\{2,3\}) &=& \pr(\tau=(2,3))+\pr(\tau=(3,2)) = \bbb{p}_{23}+\bbb{p}_{32} = \frac{1}{6}. \\
  \end{eqnarray*}
  Then, the batch-sampling $\tau$ satisfies the joint distribution $\bbb{P}$ listed in
  Table \ref{tab:jointDistribNonuniqueSampling-3-by-3}, which is still a 2nd-order tensor with constant length 3 of each dimension.
  The marginal distributions of $\tau_1$ and $\tau_2$ are listed at the last column and row of
  Table \ref{tab:jointDistribNonuniqueSampling-3-by-3}, respectively.
\renewcommand\arraystretch{1.5}
\begin{table}[htbp]
\setlength{\abovecaptionskip}{0pt}
\setlength{\belowcaptionskip}{10pt} \centering{
\caption{\label{tab:jointDistribNonuniqueSampling-3-by-3}
The joint distribution $\bbb{P}$ of non-unique batch-sampling $\tau=(\tau_1,\tau_2)$,
and the marginal distributions of $\tau_1$ and $\tau_2$: $A\in\reals^{3\times n}$, $\mathcal{G}=\{\{1\},\{2\},\{3\},\{1,2\},\{1,3\},\{2,3\}\}$.}
\begin{tabular}{|c|ccc|c|}\hline
\diagbox{$\tau_1$}{$\tau_2$} & 1 & 2 & 3 & $\tau_1=i$ \\ \hline
                           1 & $\frac{1}{6}$ & $\bbb{p}_{12}$ & $\bbb{p}_{13}$ & $\frac{1}{6}+\bbb{p}_{12}+\bbb{p}_{13}$ \\
                           2 & $\frac{1}{6}-\bbb{p}_{12}$ & $\frac{1}{6}$ & $\bbb{p}_{23}$ & $\frac{1}{3}-\bbb{p}_{12}+\bbb{p}_{23}$ \\
                           3 & $\frac{1}{6}-\bbb{p}_{13}$ & $\frac{1}{6}-\bbb{p}_{23}$ & $\frac{1}{6}$ & $\frac{1}{2}-\bbb{p}_{13}-\bbb{p}_{23}$ \\ \hline
                  $\tau_2=j$ & $\frac{1}{2}-\bbb{p}_{12}-\bbb{p}_{13}$ & $\frac{1}{3}+\bbb{p}_{12}-\bbb{p}_{23}$ & $\frac{1}{6}+\bbb{p}_{13}+\bbb{p}_{23}$ & 1 \\
\hline
\end{tabular}}
\end{table}
  According to Definition \ref{def:batchSampling}, the diagonal matrices $P_1$ and $P_2$ for $\tau_1$ and
  $\tau_2$ are of the forms
  \begin{eqnarray*}
    P_1 = \begin{bmatrix}
              \frac{1}{6}+\bbb{p}_{12}+\bbb{p}_{13} &   &   \\
                & \frac{1}{3}-\bbb{p}_{12}+\bbb{p}_{23} &   \\
                &   & \frac{1}{2}-\bbb{p}_{13}-\bbb{p}_{23} \\
            \end{bmatrix}
  \end{eqnarray*}
  and
  \begin{eqnarray*}
    P_2 = \begin{bmatrix}
              \frac{1}{2}-\bbb{p}_{12}-\bbb{p}_{13} &   &   \\
                & \frac{1}{3}+\bbb{p}_{12}-\bbb{p}_{23} &   \\
                &   & \frac{1}{6}+\bbb{p}_{13}+\bbb{p}_{23} \\
            \end{bmatrix}.
  \end{eqnarray*}
  Obviously, the property $P_1+P_2\succ 0$ gets satisfied too.

\end{example}

\end{document}